\documentclass[preprints,article,accept,moreauthors,pdftex,a4paper]{mdpi}

\firstpage{1}
\pubvolume{1}
\issuenum{1}
\articlenumber{0}
\pubyear{2022}
\copyrightyear{2022}
\datereceived{}
\dateaccepted{}
\datepublished{}
\hreflink{} % If needed use \linebreak
\pdfoutput=1

\usepackage{amssymb}
\usepackage{tabu}
\usepackage{subcaption}

\Title{Can DtN and GenEO coarse spaces be sufficiently robust for heterogeneous Helmholtz problems?}

\TitleCitation{Can DtN and GenEO coarse spaces be sufficiently robust for heterogeneous Helmholtz problems?\color{white}}

\Author{Niall Bootland $^{1,}$*\orcidA{} and Victorita Dolean $^{1,2}$\orcidB{}}

\AuthorNames{Niall Bootland and Victorita Dolean}

\AuthorCitation{Bootland, N.; Dolean, V.}
\address{%
$^{1}$ \quad Department of Mathematics and Statistics, University of Strathclyde, Glasgow, G1 1XH, UK\\
$^{2}$ \quad Laboratoire J.A.~Dieudonn\'e, CNRS, University C\^ote d'Azur, Nice, 06000, France; work@victoritadolean.com}

\corres{Correspondence: niall.bootland@strath.ac.uk}

\abstract{Numerical solution of heterogeneous Helmholtz problems presents various computational challenges, with descriptive theory remaining out of reach for many popular approaches. Robustness and scalability are key for practical and reliable solvers in large-scale applications, especially for large wave number problems. In this work we explore the use of a GenEO-type coarse space to build a two-level additive Schwarz method applicable to highly indefinite Helmholtz problems. Through a range of numerical tests on a 2D model problem, discretised by finite elements on pollution-free meshes, we observe robust convergence, iteration counts that do not increase with the wave number, and good scalability of our approach. We further provide results showing a favourable comparison with the DtN coarse space. Our numerical study shows promise that our solver methodology can be effective for challenging heterogeneous applications.}


\begin{document}
%%%%%%%%%%%%%%%%%%%%%%%%%%%%%%%%%%%%%%%%%%

\section{Introduction}
\label{sec:Introduction}

Consider solving the heterogeneous Helmholtz problem: for a bounded domain $\Omega \subset \mathbb{R}^{d}$, $d = 2, 3$, we wish to find $u(\boldsymbol{x}) \colon \Omega \rightarrow \mathbb{C}$ such that
\begin{subequations}
\label{HelmholtzSystem}
\begin{align}
	\label{HelmholtzEquation}
	-\Delta u - k^{2} u & = f & & \text{in } \Omega,\\
	\label{HelmholtzBCs}
	\mathcal{C}(u) & = 0 & & \text{on } \partial\Omega,
\end{align}
\end{subequations}
where $\mathcal{C}$ incorporates some appropriate boundary conditions. For the heterogeneous problem we suppose the wave number $k(\boldsymbol{x}) > 0$ is a function of space, defined by the ratio $k = \omega/c$ of the angular frequency $\omega$ and the wave speed $c(\boldsymbol{x})$. In this work we investigate the use of an overlapping Schwarz preconditioner with a suitably chosen coarse space based on solving local eigenvalue problems on each subdomain.

The ability to compute solutions to the Helmholtz problem \eqref{HelmholtzSystem} is important across many disciplines of science and engineering. As the prototypical model for frequency-domain wave propagation, it features within the fields of optics, acoustics and seismology amongst others. Further applications can be found in imaging science, such as through medical imaging techniques and geophysical studies of the Earth's subsurface within, for instance, the oil industry. Nonetheless, it is challenging to develop efficient computational methods to solve \eqref{HelmholtzSystem}, particularly when the wave number $k$ becomes large.

Discretisation of \eqref{HelmholtzSystem} by standard approaches, such as Lagrange finite elements as we shall use here, results in large linear systems to be solved which are indefinite, non-self-adjoint, and ill-conditioned; see also \cite{Moiola:2014:SIREV}. These systems present various difficulties to solve, especially in the presence of complex heterogeneities, at high frequencies (large $k$), or when solutions include many wavelengths in the domain. Classical methods for solving such large systems typically fail for several reasons, as detailed in \cite{Ernst:2012:NAM, Gander:2019:SIREV}, and specialist approaches must be employed for a robust solver. While much progress has been made for symmetric positive definite problems, such techniques cannot be applied out-of-the-box and extensions to tackle indefinite and non-self-adjoint problems may not be clear. This has led to a number of approaches being developed in recent years, aiming to bridge this gap. For the Helmholtz problem this includes parallel direct solvers, such as \cite{Wang:2011:OMO,Gillman:2015:BIT}, and preconditioned iterative methods that utilise multigrid, such as \cite{Calandra:2013:AIT,Hu:2016:SPD}, or the so-called ``shifted Laplace'' approach \cite{Erlangga:2004:OAC,Erlangga:2006:ANM}. This latter (complex-)shifted Laplace preconditioner has seen much interest into its practical use \cite{Cocquet:2017:HLA} and further developments through deflation techniques, most recently in \cite{Lahaye:2017:HTC,Dwarka:2020:HLA}.

Another broad class of solvers are domain decomposition methods, which provide a natural balance between using direct and iterative solvers. Specialist methods are again required for the Helmholtz problem, a popular set of which fall under the heading of ``sweeping'' methods \cite{Engquist:2011:SP1,Engquist:2011:SP2}. These are multiplicative domain decomposition methods and linked to a variety of other approaches, including optimised Schwarz methods, as detailed in the recent survey \cite{Gander:2019:SIREV}. While sweeping is a conceptually serial approach, much work has been done to incorporate parallelism. Of particular note is the ``L-Sweeps'' method \cite{Taus:2020:LSweeps}, stated to be the first parallel scalable preconditioner for high-frequency Helmholtz problems; a review of developments in this area is provided in the introductions of \cite{Taus:2020:LSweeps,Dai:2022:MSP}. Another popular approach are FETI methods, for the Helmholtz problem these include FETI-H \cite{Farhat:2000:ATL} and FETI-DPH \cite{Farhat:2005:FET}.

Within the domain decomposition community, there has also been renewed work on additive Schwarz methods, which offer a naturally parallel approach. Following on from the seminal work \cite{Despres:1990:DDM}, which utilised Robin (or impedance) transmission conditions to provide a convergent Schwarz method for the Helmholtz problem, a wealth of non-overlapping Schwarz methods have been devised; see the introduction of \cite{Claeys:2020:RTO} for a recent overview. In these methods one has to be careful to either avoid or treat cross points (where three or more subdomains meet), as can be done in the robust treatment of \cite{Claeys:2020:RTO}. Many optimised approaches rely on deriving higher-order transmission conditions, such as through second order impedance operators in \cite{Gander:2002:OSM}, absorbing boundary conditions (ABCs) \cite{Boubendir:2012:AQO}, or non-local operators \cite{Collino:2000:DDM}. Ideally one would use the Dirichlet-to-Neumann (DtN) map (a Poincar{\'e}--Steklov operator) to provide transparent transmission conditions, but this is prohibitive in practice and so these optimised Schwarz methods in essence try to approximate this operator.

Overlapping Schwarz methods for the Helmholtz problem, see for example \cite{Cai:1998:OSA,Kimn:2007:ROB,Gander:2016:OSM}, have also received renewed attention in recent years and it is this type of method we shall consider. A successful approach is to design additive Schwarz methods based on including absorption (a complex shift $k^{2} \mapsto k^{2} + i \varepsilon$), with absorption parameter $\varepsilon$; see \cite{Kimn:2013:SLR,Graham:2017:RRO}. Theoretical work to understand the effectiveness of this approach can be found in \cite{Graham:2017:DDP,Graham:2020:DDW} and for the heterogeneous problem in \cite{Gong:2021:DDP}; see also \cite{Gong:2021:COP}. To be scalable, such additive Schwarz methods require a second level, known as a coarse space (though see \cite{Bootland:2022:AOP} for a novel analysis for the absorptive problem). For the two-level methods in \cite{Graham:2017:DDP,Graham:2017:RRO} this is provided through a coarse grid; an approach which is effective also for the time-harmonic Maxwell problem \cite{Bonazzoli:2019:AOP}.

As well as to provide scalability, coarse spaces have been devised to provide robustness to heterogeneity. This is exemplified by the ``Generalised Eigenproblems in the Overlap'' (GenEO) approach for symmetric positive definite (SPD) problems \cite{Spillane:2014:ARC}. This approach provides a spectral coarse space, where appropriate local eigenvalue problems are solved to provide a two-level method. Another spectral coarse space is the DtN coarse space \cite{Nataf:2011:ACS}, which has been extended and investigated for the Helmholtz problem in \cite{Conen:2014:ACS,Bootland:2019:ODN}. While the standard GenEO theory applies only in the SPD case (though see \cite{Bootland:2021:OSM}), in this work we develop and explore a GenEO-type method for the Helmholtz problem and show numerically that, for a 2D model problem of a wave guide, it provides a scalable approach that is robust to heterogeneity and increasing wave number in terms of the iteration count of a preconditioned GMRES method. Companion results for large benchmark problems on coarsely resolved meshes that typically arise in applications, along with comparisons to other methods, are found in \cite{Bootland:2021:ACO}. Theoretical results on variants of DtN and GenEO for Helmholtz problems are currently out of reach but promising numerical results, based on heuristics, were obtained in \cite{Conen:2014:ACS,Bootland:2019:ODN}, showing the potential of these methods in practice.

The primary aim of this work is to explore the utility of a GenEO-type method for the heterogeneous Helmholtz problem \eqref{HelmholtzSystem}; we call this approach H-GenEO. In particular, we highlight the following contributions:
\begin{itemize}
	\item We present a range of numerical tests, on pollution-free meshes, comparing our proposed H-GenEO approach with another spectral coarse space applicable to the Helmholtz problem, namely the DtN method.
	\item We investigate the use of appropriate thresholding for the required generalised eigenproblems in both the DtN and H-GenEO coarse spaces.
	\item We consider robustness to non-uniform decomposition, heterogeneity, and increasing wave number as well as the scalability of the methods. We find that only the H-GenEO approach is scalable and robust to all of these factors for a 2D model problem.
	\item We provide both weak and strong scalability tests for H-GenEO applied to high wave number problems.
\end{itemize}

The remainder of this work is structured as follows. We begin by considering a finite element discretisation of the Helmholtz problem in Section~\ref{sec-sub:FEDiscretisation}, before outlining the underlying domain decomposition methodology we use in Section~\ref{sec-sub:DD}. The main topic of interest, that of suitable spectral coarse spaces for the heterogeneous Helmholtz problem, is then detailed in Section~\ref{sec-sub:CoarseSpaces}. Extensive numerical results on a 2D model problem are provided in Section~\ref{sec:Results}, along with a discussion of our findings. Finally, we draw together our conclusions in Section~\ref{sec:Conclusions}.

%%%%%%%%%%%%%%%%%%%%%%%%%%%%%%%%%%%%%%%%%%
\section{Materials and Methods}
\label{sec:Methods}

\subsection{Finite Element Discretisation}
\label{sec-sub:FEDiscretisation}

The problem we consider in this work is the interior Helmholtz problem \eqref{HelmholtzSystem}, for which we must prescribe appropriate boundary conditions. In practical applications, the computational domain $\Omega$ is often truncated and the physically relevant condition, namely the far field Sommerfeld radiation condition, must be approximated on the non-physical boundary of $\Omega$. This allows for appropriate wave behaviour to be modelled in a bounded domain. This simplest approximation which is widely used is that of a Robin (or impedance) condition and this is what we shall consider; other approaches include ABCs \cite{Zarmi:2013:AGA} or perfectly matched layers (PML) \cite{Harari:2000:AAN,Beriot:2021:AAP}. We also suppose that Dirichlet conditions may be imposed on a boundary $\Gamma_{D} \subset \partial\Omega$, with the Robin condition on the remaining boundary $\Gamma_{R} = \partial\Omega \setminus \Gamma_{D}$. Thus, in general, we seek the solution of the boundary value problem
\begin{subequations}
\label{HelmholtzSpecificSystem}
\begin{align}
	\label{HelmholtzSpecificEquation}
	-\Delta u - k^{2} u & = f & & \text{in } \Omega,\\
	\label{HelmholtzSpecificDirichletBC}
	u & = u_{\Gamma_{D}} & & \text{on } \Gamma_{D},\\
	\label{HelmholtzSpecificRobinBC}
	\frac{\partial u}{\partial \boldsymbol{n}} + i k u & = 0 & & \text{on } \Gamma_{R},
\end{align}
\end{subequations}
where the forcing function $f(\boldsymbol{x})$ incorporates any sources in the domain. Note that if $\Gamma_{R} \neq \emptyset$ the problem is well posed but if $\Gamma_{R} = \emptyset$ the problem is ill-posed for certain choices of $k$ related to eigenfunctions of the Laplacian.

To discretise \eqref{HelmholtzSpecificSystem} we use standard Lagrange finite elements; the details can be found in, for example, \cite{Bootland:2021:ACO} and so we provide only an outline here. Defining the relevant trial and test spaces $V = \left\lbrace u \in H^{1}(\Omega) \colon u = u_{\Gamma_{D}}\text{ on } \Gamma_{D}\right\rbrace$ and $V_{0} = \left\lbrace u \in H^{1}(\Omega) \colon u = 0\text{ on } \Gamma_{D}\right\rbrace$, the weak formulation of \eqref{HelmholtzSpecificSystem} is to find $u \in V$ such that
\begin{align}
\label{WeakForm}
a(u,v) & = F(v) & & \forall \ v \in V_{0},
\end{align}
where
\begin{align}
\label{WeakFormTerms}
a(u,v) & = \int_{\Omega} \left( \nabla u \cdot \nabla \bar{v} - k^2 u \bar{v}\right) \, \mathrm{d}\boldsymbol{x} + \int_{\Gamma_{R}} i k u \bar{v} \, \mathrm{d}s & & \text{and} & F(v) & = \int_{\Omega} f \bar{v} \, \mathrm{d}\boldsymbol{x}.
\end{align}
Assuming a simplicial mesh $\mathcal{T}^{h}$ of $\Omega$ with characteristic element diameter $h$, piecewise polynomial finite element approximation reduces the problem to solving the complex linear system
\begin{align}
\label{LinearSystem}
A\mathbf{u} = \mathbf{f},
\end{align}
with coefficient matrix $A \in \mathbb{C}^{n \times n}$ and right-hand side vector $\mathbf{f} \in \mathbb{C}^{n}$ stemming from $a(\cdot,\cdot)$ and $F(\cdot)$ respectively; see \cite{Bootland:2021:ACO}. For accurate discretisation of the Helmholtz problem, the number of degrees of freedom $n$ is required to be large, especially for large $k$. Indeed, to maintain the same level of accuracy of discrete solutions as $k$ increases then the number of mesh points must increase faster than $k$, due to the pollution effect \cite{Babuska:1997:IPE}. This depends on the polynomial order of the approximation used: for instance, using piecewise linear (P1) finite elements $k^{3} h^{2}$ must be bounded and so $h$ must shrink as $\mathcal{O}(k^{-3/2})$. Higher order finite elements can reduces this restriction on $h$ but ultimately the interpolation properties of such methods degrade. Here we will utilise standard P1 elements and maintain a discretisation such that $k^{3} h^{2}$ is fixed to avoid the pollution effect.

\subsection{Underlying Domain Decomposition Method}
\label{sec-sub:DD}

To solve the discrete Helmholtz problem \eqref{LinearSystem} we will use GMRES accelerated via a two-level overlapping domain decomposition preconditioner. For the underlying one-level method we consider the optimised restricted additive Schwarz (ORAS) method \cite{Dolean:15:DDM} (sometimes known as WRAS-H \cite{Kimn:2007:ROB}; see discussions in \cite{Gong:2021:DDP,Gong:2021:COP} and \cite{Graham:2017:DDP,Graham:2017:RRO} for a related IMPHRAS1 method). They key difference in ORAS compared with standard additive Schwarz methods such as RAS is that the local Dirichlet problems on subdomains are replaced by appropriate Robin problems.

To formulate the domain decomposition preconditioner, we first suppose that $\Omega$ is decomposed into non-overlapping subdomains $\left\lbrace\Omega'_{s}\right\rbrace_{s=1}^{N}$, assumed to be resolved by the mesh $\mathcal{T}^{h}$. To give an overlapping decomposition, a layer of adjoining mesh elements is added to give overlapping subdomains $\left\lbrace\Omega_{s}\right\rbrace_{s=1}^{N}$ by way of the extension
\begin{align}
\label{OverlappingSubdomains}
\Omega_{s} = \mathrm{Int}\left(\bigcup_{\mathrm{supp}(\phi_{j})\cap\Omega'_{s}\neq\emptyset} \mathrm{supp}(\phi_{j})\right),
\end{align}
where $\left\lbrace\phi_{j}\right\rbrace_{j=1}^{n}$ are the nodal basis functions of the finite element space, $\mathrm{Int}(\cdot)$ denotes the interior of a domain and $\mathrm{supp}(\cdot)$ the support of a function. Further layers of elements can be added in a recursive manner in order to obtain subdomains with larger overlap, if desired.

Once we have an overlapping decomposition into subdomains $\left\lbrace\Omega_{s}\right\rbrace_{s=1}^{N}$, we define the required operators for the Schwarz preconditioner. We let $R_{s} \in \mathbb{R}^{n_{s} \times n}$ be the discrete form of the restriction operator, restricting functions to the subdomain $\Omega_{s}$, where $n_{s}$ is the number of degrees of freedom in $\Omega_{s}$. The corresponding extension operator, $R_{s}^{T}$, then acts as an extension by zero outside of $\Omega_{s}$. We also utilise a partition of unity, having the discrete form of a diagonal matrix $D_{s} \in \mathbb{R}^{n_{s} \times n_{s}}$ satisfying $\sum_{s=1}^{N} R_{s}^{T} D_{s} R_{s} = I$, which appropriately scales the multiple subdomain contributions in the overlapping regions. Finally, within ORAS we require the solution of local Robin problems given by
\begin{subequations}
\label{ORASLocalSystem}
\begin{align}
	\label{ORASLocalHelmholtzEquation}
	-\Delta w_{s} - k^{2} w_{s} & = f & & \text{in } \Omega_{s},\\
	\label{ORASLocalRobinBC}
	\frac{\partial w_{s}}{\partial \boldsymbol{n}_{s}} + i k w_{s} & = 0 & & \text{on } \partial\Omega_{s}\setminus\partial\Omega,\\
	\label{ORASLocalProblemBCs}
	\mathcal{C}(w_{s}) & = 0 & & \text{on } \partial\Omega_{s}\cap\partial\Omega,
\end{align}
\end{subequations}
with $\mathcal{C}$ representing the underlying problem boundary conditions on $\partial\Omega$, namely \eqref{HelmholtzSpecificDirichletBC}--\eqref{HelmholtzSpecificRobinBC}. Note the use of the Robin condition \eqref{ORASLocalRobinBC} ensures the solvability of these local problems. Defining the equivalent finite element discretisation of \eqref{ORASLocalSystem} to be given by the stiffness matrix $\widehat{A}_{s} \in \mathbb{R}^{n_{s} \times n_{s}}$, the construction of the one-level ORAS preconditioner is given by the sum
\begin{align}
\label{ORAS}
M_{\text{ORAS}}^{-1} = \sum_{s=1}^{N} R_{s}^{T} D_{s} \widehat{A}_{s}^{-1} R_{s}.
\end{align}
Note that the local solutions on each subdomain given by $\widehat{A}_{s}^{-1}$ can be carried out in parallel.

In order to provide robustness and scalability, the ORAS method above must be augmented by use of a coarse space to provide a two-level method. A coarse space can be thought of as a collection of linearly independent column vectors $Z$. The vectors that are incorporated into $Z$ are key to providing scalability, especially for indefinite problems, such as the Helmholtz problems we solve here, where the addition of a coarse space need not improve performance of the underlying one-level method \cite{Fish:2000:GBT}. There are several ways to incorporate the coarse space; here we consider an effective approach that is based on deflation. For this, a coarse space operator $E = Z^{\dagger} A Z$ is constructed as well as the coarse correction operator $Q = Z E^{-1} Z^{\dagger}$, which we incorporate to give a two-level ORAS method
\begin{align}
\label{2LevelAdaptiveDeflationPreconditioner}
M_{\text{ORAS,2}}^{-1} = M_{\text{ORAS}}^{-1}(I-AQ) + Q.
\end{align}
We now turn our attention the choice of coarse space.

\subsection{Spectral Coarse Spaces}
\label{sec-sub:CoarseSpaces}

In this work we consider and explore spectral coarse spaces for the discrete Helmholtz problem \eqref{LinearSystem}. These utilise local eigenvalue problems on subdomains in order to build a global coarse space. We will review the DtN coarse space \cite{Conen:2014:ACS} before detailing a new coarse space for Helmholtz problems based on GenEO technology \cite{Spillane:2014:ARC}. We will then show a link between these two approaches.

\begin{Remark}[Notation]
We utilise the following notation for local Dirichlet, Robin, and Neumann matrices: For a variational problem which gives rise to a system matrix $B$, we denote by $B_{s}$ the corresponding local Dirichlet matrix on $\Omega_{s}$. In the case that Robin conditions are used on internal subdomain interfaces the local problem matrix is denoted by $\widehat{B}_{s}$. On the other hand, if Neumann conditions are used on such interfaces we denote the local matrix by $\widetilde{B}_{s}$.
\end{Remark}

\subsubsection{The DtN Coarse Space}

The Dirichlet-to-Neumann (DtN) coarse space, first studied in \cite{Nataf:2010:ATL,Nataf:2011:ACS} for elliptic problems, is based on solving local eigenvalue problems on subdomain boundaries related to a DtN map. Harmonic extensions of low-frequency modes on subdomains are then used to provide a coarse space. In order to define this approach for the Helmholtz problem, as explored in \cite{Conen:2014:ACS}, we first require the Helmholtz extension operator from the subdomain boundary $\partial\Omega_{s}$.

On each subdomain, let $\Gamma_{s} = \partial\Omega_{s} \setminus \partial\Omega$ and suppose we have Dirichlet data $v_{\Gamma_{s}}$ on $\Gamma_{s}$. The Helmholtz extension $v$ into $\Omega_{s}$ is given by solving
\begin{subequations}
\label{HelmholtzExtension}
\begin{align}
\label{HelmholtzExtensionEquation}
-\Delta v - k^{2} v & = 0 & & \text{in } \Omega_{s},\\
\label{HelmholtzExtensionDirichletBC}
v & = v_{\Gamma_{s}} & & \text{on } \Gamma_{s},\\
\label{HelmholtzExtensionProblemBC}
\mathcal{C}(v) & = 0 & & \text{on } \partial\Omega_{s}\cap\partial\Omega,
\end{align}
\end{subequations}
where $\mathcal{C}(v) = 0$ represents the original problem boundary conditions, as in \eqref{ORASLocalProblemBCs}. The DtN map takes Dirichlet data $v_{\Gamma_{s}}$ to the corresponding Neumann data on $\Gamma_{s}$, namely
\begin{align}
\label{DtNMap}
\mathrm{DtN}_{\Omega_{s}}(v_{\Gamma_{s}}) = \left.\frac{\partial v}{\partial n} \right\rvert_{\Gamma_{s}}
\end{align}
where $v$ is the Helmholtz extension defined by \eqref{HelmholtzExtension}. The associated local DtN eigenproblem on subdomain $\Omega_{s}$ is given by
\begin{align}
\label{DtNEigenproblem}
\mathrm{DtN}_{\Omega_{s}}(u_{\Gamma_{s}}) = \lambda u_{\Gamma_{s}},
\end{align}
for eigenfunctions $u_{\Gamma_{s}}$ and eigenvalues $\lambda \in \mathbb{C}$. In order to build the coarse space we take the Helmholtz extension of $u_{\Gamma_{s}}$ in $\Omega_{s}$ and extend by zero into the whole domain $\Omega$ using the partition of unity; see \cite{Conen:2014:ACS}.

To formulate the discrete version of the eigenproblems to be solved we require the coefficient matrices $\widetilde{A}_{s}$, corresponding to local Neumann problems on $\Omega_{s}$ with boundary conditions $\mathcal{C}=0$ on $\partial\Omega_{s}\cap\partial\Omega$, similar to that supplying the local Robin problems in \eqref{ORASLocalSystem}. Furthermore, we must distinguish between degrees of freedom on the boundary and the interior of the subdomain $\Omega_{s}$ and so we let $\Gamma_{s}$ and $I_{s}$ be the set of indices on the boundary and interior respectively. Recalling that  $\left\lbrace\phi_{j}\right\rbrace$ are our nodal basis functions, we also define
\begin{align}
\label{DtNMassMatrix}
M_{\Gamma_{s}} = \left(\int_{\Gamma_{s}} \phi_{j} \phi_{i} \right)_{i,j \in \Gamma_{s}}
\end{align}
to be the mass matrix on the subdomain interface. Using standard block notation to denote submatrices of $A_{s}$ and $\widetilde{A}_{s}$ the discrete DtN eigenproblem can be written as
\begin{align}
\label{DiscreteDtNEigenproblem}
\left(\widetilde{A}_{\Gamma_{s},\Gamma_{s}} - A_{\Gamma_{s},I_{s}}A_{I_{s},I_{s}}^{-1}A_{I_{s},\Gamma_{s}}\right) \mathbf{u}_{\Gamma_{s}} = \lambda M_{\Gamma_{s}} \mathbf{u}_{\Gamma_{s}}.
\end{align}
We then make use of the Helmholtz extension of $\mathbf{u}_{\Gamma_{s}}$ to degrees of freedom in $I_{s}$ given by $\mathbf{u}_{I_{s}} = - A_{I_{s},I_{s}}^{-1}A_{I_{s},\Gamma_{s}} \mathbf{u}_{\Gamma_{s}}$. Letting $\mathbf{u}_{s}$ denote the complete local vector representing the Helmholtz extension, the corresponding global vector which enters the coarse space Z is $R_{s}^{T} D_{s} \mathbf{u}_{s}$. Further motivation and details on the DtN eigenproblems can be found in \cite{Conen:2014:ACS}.

What remains is to determine which eigenvectors of \eqref{DiscreteDtNEigenproblem} should be incorporated into the coarse space. A variety of selection criteria were investigated in \cite{Conen:2014:ACS} which made it clear that the best choice was to select eigenvectors corresponding to eigenvalues with the smallest real part. That is, a threshold on the abscissa $\eta = \mathrm{Re}(\lambda)$ should be used, namely
\begin{align}
\label{Threshold}
\eta < \eta_{\text{max}},
\end{align}
where $\eta_{\text{max}}$ depends on $k_{s} = \max_{\vec{x}\in\Omega_{s}} k(\vec{x})$. The choice $\eta_{\text{max}} = k_{s}$ is advocated in \cite{Conen:2014:ACS}, however, we recently showed that taking a slightly larger threshold $\eta_{\text{max}} = k_{s}^{4/3}$ can be beneficial in certain cases in order to gain robustness to the wave number \cite{Bootland:2019:ODN}. Unfortunately, this only occurs for the homogeneous problem with sufficiently uniform subdomains. To construct a more robust coarse space we will build upon the GenEO approach.

\subsubsection{The GenEO Coarse Space}

The Generalised Eigenproblems in the Overlap (GenEO) coarse space was derived in \cite{Spillane:2014:ARC} to provide a rigorously robust approach for symmetric positive definite problems even in the presence of heterogeneities. In recent years, this approach have been extended and used within various settings and applications, for example \cite{Dolean:15:DDM,Haferssas:2017:AAS,Nataf:2019:AGD,Bootland:2021:OSM}; see also the discussion on developments for other spectral coarse spaces in \cite{Spillane:2021:AAT}.

Within the original derivation \cite{Spillane:2014:ARC}, the generalised eigenproblems are defined in a variational framework on the assumption that $a(\cdot,\cdot)$ is a symmetric and coercive bilinear form. On a subdomain $\Omega_{s}$, an overlapping zone $\Omega_{s}^{\circ}$ is defined as the parts of $\Omega_{s}$ which overlap with another subdomain and the local eigenproblem
\begin{align}
\label{GenEOVariationalEigenproblem}
a_{\Omega_{s}}(u,v) & = \lambda a_{\Omega_{s}^{\circ}}\left(\Xi_{s}(u),\Xi_{s}(v)\right) & & \forall \ v \in V(\Omega_{s}),
\end{align}
is solved for small $\lambda$, where $\Xi_{s}$ represents the action of the partition of unity operator on $\Omega_{s}$. To be clear, $a_{D}(\cdot,\cdot)$ represents the underlying variational problem on the domain $D$ with problem boundary conditions on $\partial\Omega$ and natural (``do nothing'') conditions on $\partial D \setminus \partial\Omega$. The eigenproblem \eqref{GenEOVariationalEigenproblem} provides an appropriate link in order to bound the condition number of the preconditioned operator independently of the heterogeneity and number of subdomains. This bound depends on the smallest eigenvalue $\lambda$ whose corresponding eigenfunction is not incorporated into the GenEO coarse space. Hence, to achieve a desired rate of convergence, all eigenfunctions corresponding to eigenvalues smaller than a threshold (say, $\lambda < \lambda_{\text{max}}$) must be computed. Alternatively, one may opt in practice to compute a fixed number of eigenfunctions per subdomain and use these in order to accelerate convergence.

The restriction within the right-hand side of \eqref{GenEOVariationalEigenproblem} to the overlapping zone is not an essential requirement and alternative formulations can be used. In order to remove the need to track overlap regions, one possibility is to replace $\Omega_{s}^{\circ}$ with the whole of $\Omega_{s}$, as in \cite{Dolean:15:DDM,Haferssas:2017:AAS}. When formulating the discrete eigenproblem, this requires only the local Neumann matrix $\widetilde{A}_{s}$ to be constructed. In this case we must solve
\begin{align}
\label{DiscreteGenEOEigenproblem}
\widetilde{A}_{s} \mathbf{u} = \lambda D_{s} A_{s} D_{s} \mathbf{u},
\end{align}
where $A_{s} = R_{s} A R_{s}^{T}$ is the local Dirichlet matrix, a sub-matrix of $A$. As with the DtN method, the vectors which then go into the coarse space $Z$ are $R_{s}^{T} D_{s} \mathbf{u}$. It is this form of the GenEO eigenproblem in \eqref{DiscreteGenEOEigenproblem} that we shall build upon to develop an approach tailored to the heterogeneous Helmholtz problem.

\subsubsection{H-GenEO: a GenEO-Type Coarse Space for Helmholtz Problems}

In consideration of GenEO approaches for the Helmholtz problem, a key hurdle is the loss of operators being definite or self-adjoint. While some progress has been made at the theoretical level for closely related problems \cite{Bootland:2022:GCS,Bootland:2021:OSM}, currently available analysis has yet to overcome all the challenges present for the Helmholtz problem and rigorous justification for a choice spectral coarse space remains out of reach. The approach taken in \cite{Bootland:2022:GCS,Bootland:2021:OSM} is to formulate the GenEO eigenproblem for a nearby symmetric positive definite problem, here corresponding to a Laplace problem and hence given the name $\Delta$-GenEO. This problem can be given by setting $k=0$ in \eqref{HelmholtzSpecificSystem} and \eqref{WeakFormTerms}. Letting $L_{s}$ be the local Dirichlet matrix for this problem in $\Omega_{s}$, and $\widetilde{L}_{s}$ the equivalent Neumann matrix, the $\Delta$-GenEO eigenproblem is given by
\begin{align}
\label{DiscreteLaplaceGenEOEigenproblem}
\widetilde{L}_{s} \mathbf{u} = \lambda D_{s} L_{s} D_{s} \mathbf{u}.
\end{align}
This is a positive (semi-)definite eigenproblem with real non-negative eigenvalues $\lambda$ and, as such, eigenvectors can be chosen in the standard way using a threshold $\lambda < \lambda_{\text{max}}$. Unfortunately, for the Helmholtz problem the $\Delta$-GenEO approach can perform rather poorly when $k$ becomes large (see Table~\ref{Table:WideMethodsN25} in Section~\ref{sec:Results}), as might be anticipated from the fact that solutions to the Laplace problem differ considerably to those of the Helmholtz problem in this range.

In order to provide an appropriate spectral coarse space for the Helmholtz problem it stands to reason that a Helmholtz operator must be included. If we try to apply the GenEO eigenproblem \eqref{DiscreteGenEOEigenproblem} as is, with matrices stemming from the Helmholtz bilinear form in \eqref{WeakFormTerms}, we must first note that the problem is non-self-adjoint and, as such, eigenvalues $\lambda$ are no longer real in general. As a threshold criterion, as with the DtN approach, we can consider the abscissa $\eta = \mathrm{Re}(\lambda)$ instead and seek eigenvectors corresponding to $\eta < \eta_{\text{max}}$. Unfortunately, this formulation fails to be robust and numerically we have found that the eigensolver we use can often break down. Given that, without any relevant theory, it is no longer clear \eqref{DiscreteGenEOEigenproblem} provides an appropriate eigenproblem, we develop a different approach which yields a robust method in our numerical experiments.

While the above approaches consider only Helmholtz, or only Laplace, operators in their formulation, we instead link the underlying Helmholtz problem to the positive definite Laplace problem. Since this GenEO-type method targets the Helmholtz problem we call it ``H-GenEO''. The local eigenproblem utilised is given by
\begin{align}
\label{DiscreteH-GenEOEigenproblem}
\widetilde{A}_{s} \mathbf{u} = \lambda D_{s} L_{s} D_{s} \mathbf{u}.
\end{align}
Since eigenvalues are complex (though, in our experience, tend to cluster close to the real line) we threshold based on the abscissa $\eta < \eta_{\text{max}}$. We will see the effectiveness of the H-GenEO coarse space for a 2D model problem in Section~\ref{sec:Results}. When used with a threshold (typically we use $\eta_{\text{max}} = \frac{1}{2}$), we observe that the number of GMRES iterations for the H-GenEO approach is independent of the wave number $k$, which is not true of the DtN method in general. Before continuing to our numerical results, we first show a link between the DtN and GenEO eigenproblems and how the H-GenEO method is related.

\subsubsection{A Link Between DtN and GenEO}

While the DtN and GenEO eigenproblems look rather different, here we show a link between the two approaches. In particular, consider the GenEO method when we remove the partition of unity matrices $D_{s}$ from \eqref{DiscreteGenEOEigenproblem}, so that we are left with the Neumann matrix on the left and the Dirichlet matrix on the right, namely $\widetilde{A}_{s}\mathbf{u} = \lambda A_{s} \mathbf{u}$. Further, let us use subscripts $O_{s}$ and $I_{s}$ to denote overlap and interior degrees of freedom in $\Omega_{s}$ respectively. Then, moving all terms to the left-hand side, the GenEO eigenproblem can be written as
\begin{align}
\label{GenEOwithoutPoU}
\left( \begin{array}{cc}
	\widetilde{A}_{O_{s},O_{s}} - \lambda A_{O_{s},O_{s}} & (1-\lambda) A_{O_{s},I_{s}} \\
	(1-\lambda) A_{I_{s},O_{s}} & (1-\lambda) A_{I_{s},I_{s}}
\end{array} \right) \left( \begin{array}{c}
	\mathbf{u}_{O_{s}} \\
	\mathbf{u}_{I_{s}}
\end{array} \right) = \left( \begin{array}{c}
	\mathbf{0} \\
	\mathbf{0}
\end{array} \right),
\end{align}
where we have used the fact that the Neumann matrix is equal to Dirichlet matrix except in the block associated with the overlap degrees of freedom. Forming the Schur complement of \eqref{GenEOwithoutPoU} with respect to the $A_{I_{s},I_{s}}$ block and dividing by $1-\lambda$ we obtain
\begin{align}
\label{GenEOwithoutPoUSchurComplement}
\left(\widetilde{A}_{O_{s},O_{s}} - A_{O_{s},I_{s}} A_{I_{s},I_{s}}^{-1} A_{I_{s},O_{s}}\right) \mathbf{u}_{O_{s}} = \frac{\lambda}{1-\lambda} \left(A_{O_{s},O_{s}} - \widetilde{A}_{O_{s},O_{s}}\right)\mathbf{u}_{O_{s}}.
\end{align}
Now if the overlap degrees of freedom $O_{s}$ are precisely those used as boundary degrees of freedom $\Gamma_{s}$ in the DtN method we see that \eqref{GenEOwithoutPoUSchurComplement} resembles the DtN eigenproblem \eqref{DiscreteDtNEigenproblem}. The primary difference stems from the right-hand side where we now have the difference between the Dirichlet and Neumann matrices on the boundary degrees of freedom as opposed to a mass matrix, though we note that these can coincide for certain choices of discretisation, for example a simple finite difference scheme on a Cartesian grid where forward or backward differences are used to approximate the Neumann condition. The other difference, aside from the partition of unity in the true GenEO approach, stems from a transformation of the eigenvalues, namely the comparative eigenvalues for the DtN method are $\mu = \frac{\lambda}{1-\lambda}$.

Now consider the H-GenEO eigenproblem \eqref{DiscreteH-GenEOEigenproblem} without the partition of unity matrices. Splitting for the moment $A = L - k^2 M$ we have
\begin{align}
\label{H-GenEOwithoutPoU}
\left( \begin{array}{cc}
	\widetilde{L}_{O_{s},O_{s}} - \lambda L_{O_{s},O_{s}} - k^2 M_{O_{s},O_{s}} & (1-\lambda) L_{O_{s},I_{s}} - k^2 M_{O_{s},I_{s}} \\
	(1-\lambda) L_{I_{s},O_{s}} - k^2 M_{I_{s},O_{s}} & (1-\lambda) L_{I_{s},I_{s}} - k^2 M_{I_{s},I_{s}}
\end{array} \right) \left( \begin{array}{c}
	\mathbf{u}_{O_{s}} \\
	\mathbf{u}_{I_{s}}
\end{array} \right) = \left( \begin{array}{c}
	\mathbf{0} \\
	\mathbf{0}
\end{array} \right).
\end{align}
Now dividing by $1-\lambda$ and defining $\kappa^2 = \frac{k^2}{1-\lambda}$ and the corresponding Helmholtz matrices with wave number $\kappa$ as $B = L - \kappa^{2}M$ we have
\begin{align}
\label{H-GenEOwithoutPoU-Rewritten}
\left( \begin{array}{cc}
	\widetilde{B}_{O_{s},O_{s}} + \frac{\lambda}{1-\lambda} \left( \widetilde{L}_{O_{s},O_{s}} - L_{O_{s},O_{s}} \right) & B_{O_{s},I_{s}} \\
	B_{I_{s},O_{s}} & B_{I_{s},I_{s}}
\end{array} \right) \left( \begin{array}{c}
	\mathbf{u}_{O_{s}} \\
	\mathbf{u}_{I_{s}}
\end{array} \right) = \left( \begin{array}{c}
	\mathbf{0} \\
	\mathbf{0}
\end{array} \right).
\end{align}
Using the fact that $\widetilde{L}_{O_{s},O_{s}} - L_{O_{s},O_{s}} = \widetilde{B}_{O_{s},O_{s}} - B_{O_{s},O_{s}}$ the resulting Schur complement system is then given as
\begin{align}
\label{H-GenEOwithoutPoUSchurComplement}
\left(\widetilde{B}_{O_{s},O_{s}} - B_{O_{s},I_{s}} B_{I_{s},I_{s}}^{-1} B_{I_{s},O_{s}}\right) \mathbf{u}_{O_{s}} = \frac{\lambda}{1-\lambda} \left(B_{O_{s},O_{s}} - \widetilde{B}_{O_{s},O_{s}}\right)\mathbf{u}_{O_{s}}.
\end{align}
Thus we see that H-GenEO corresponds to solving GenEO problems based on a wave number $\kappa$ which varies with the eigenvalue. Since we are predominantly looking for small eigenvalues $\lambda$, this wave number $\kappa = k (1-\lambda)^{-1/2}$ will then be close to $k$.

\begin{figure}[H]
\widefigure
\centering
\rotatebox{90}{\hspace*{1.65cm} DtN}
\hfill
\begin{subfigure}[b]{0.18\textwidth}
	\centering
	\includegraphics[width=\textwidth,trim=7cm 0cm 7cm 0cm ,clip]{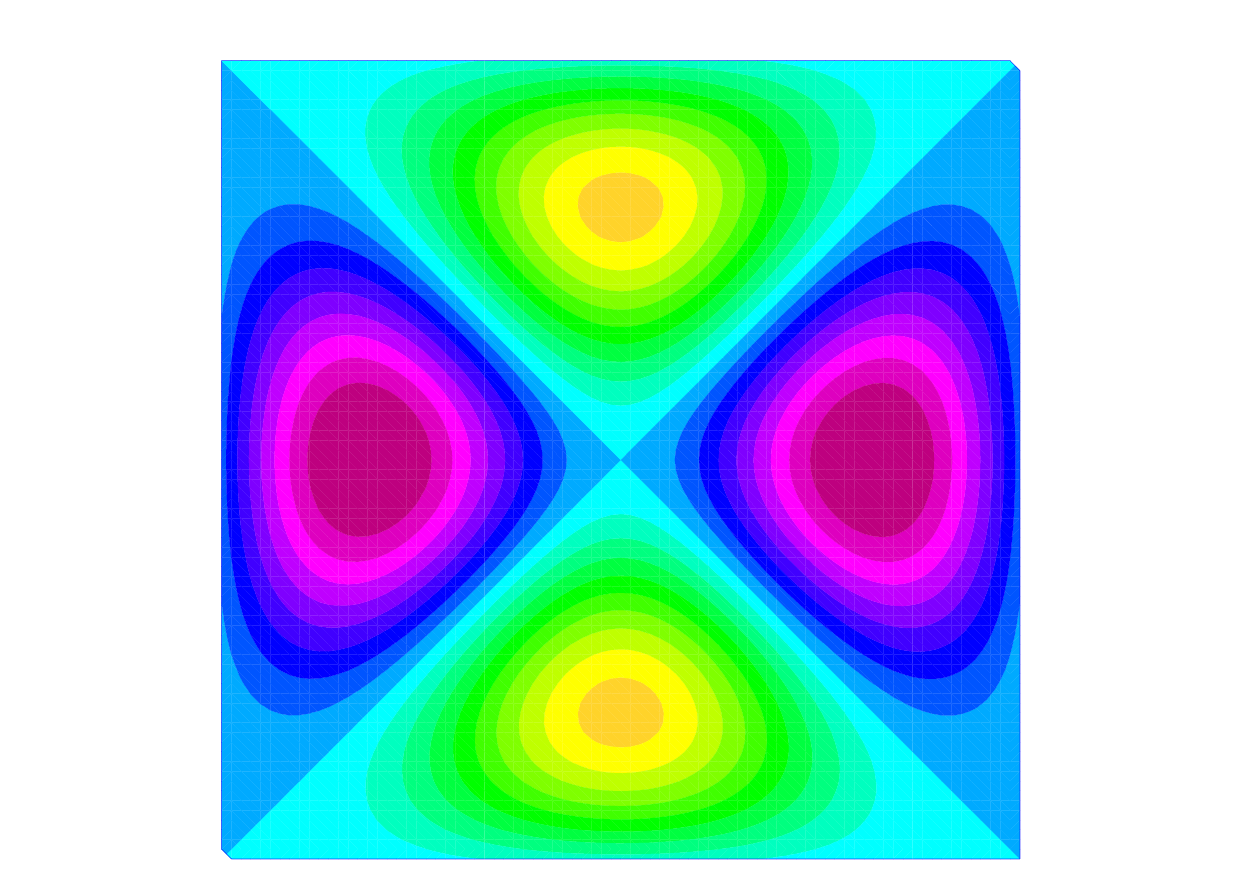}
	\caption{$\lambda_{1} = -225.2$}
\end{subfigure}
\hfill
\begin{subfigure}[b]{0.18\textwidth}
	\centering
	\includegraphics[width=\textwidth,trim=7cm 0cm 7cm 0cm ,clip]{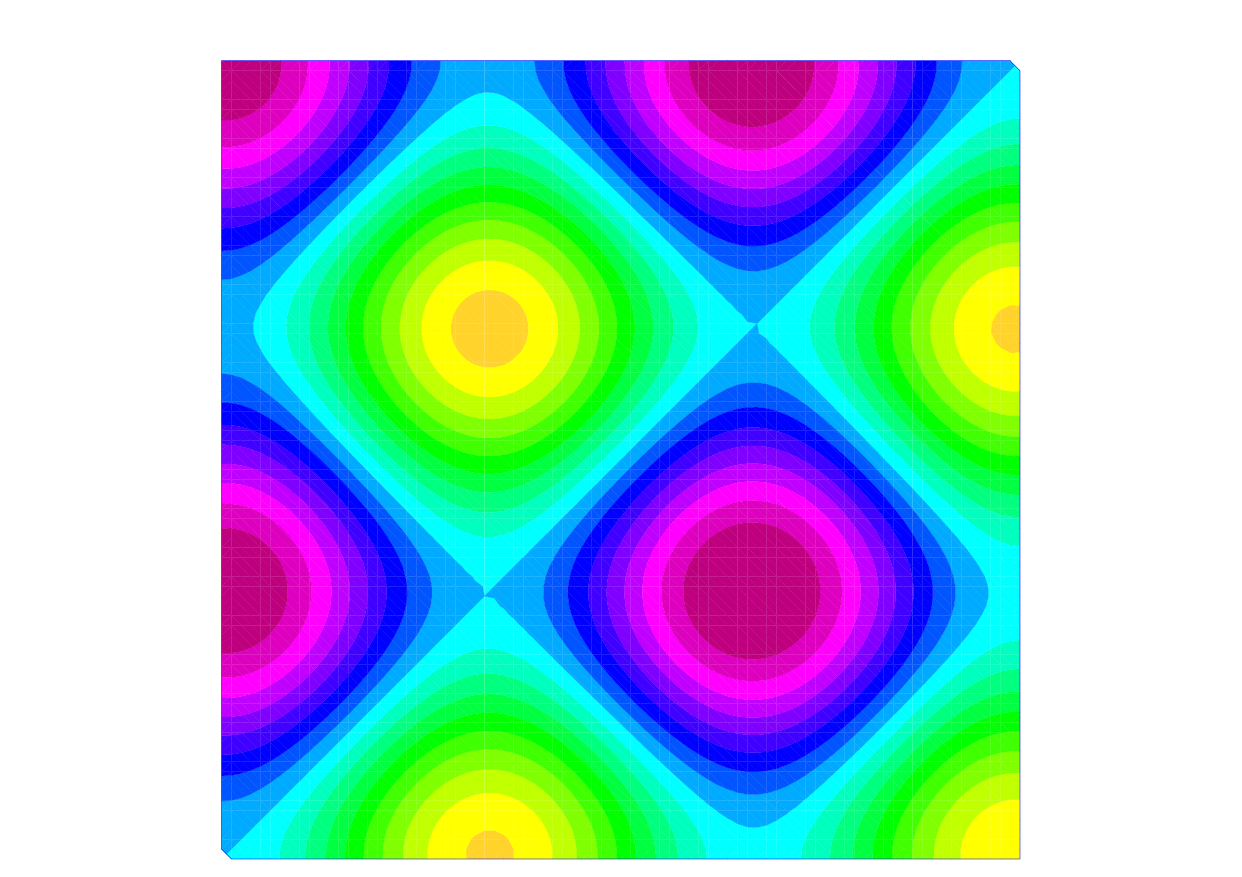}
	\caption{$\lambda_{6} = -1.6$}
\end{subfigure}
\hfill
\begin{subfigure}[b]{0.18\textwidth}
	\centering
	\includegraphics[width=\textwidth,trim=7cm 0cm 7cm 0cm ,clip]{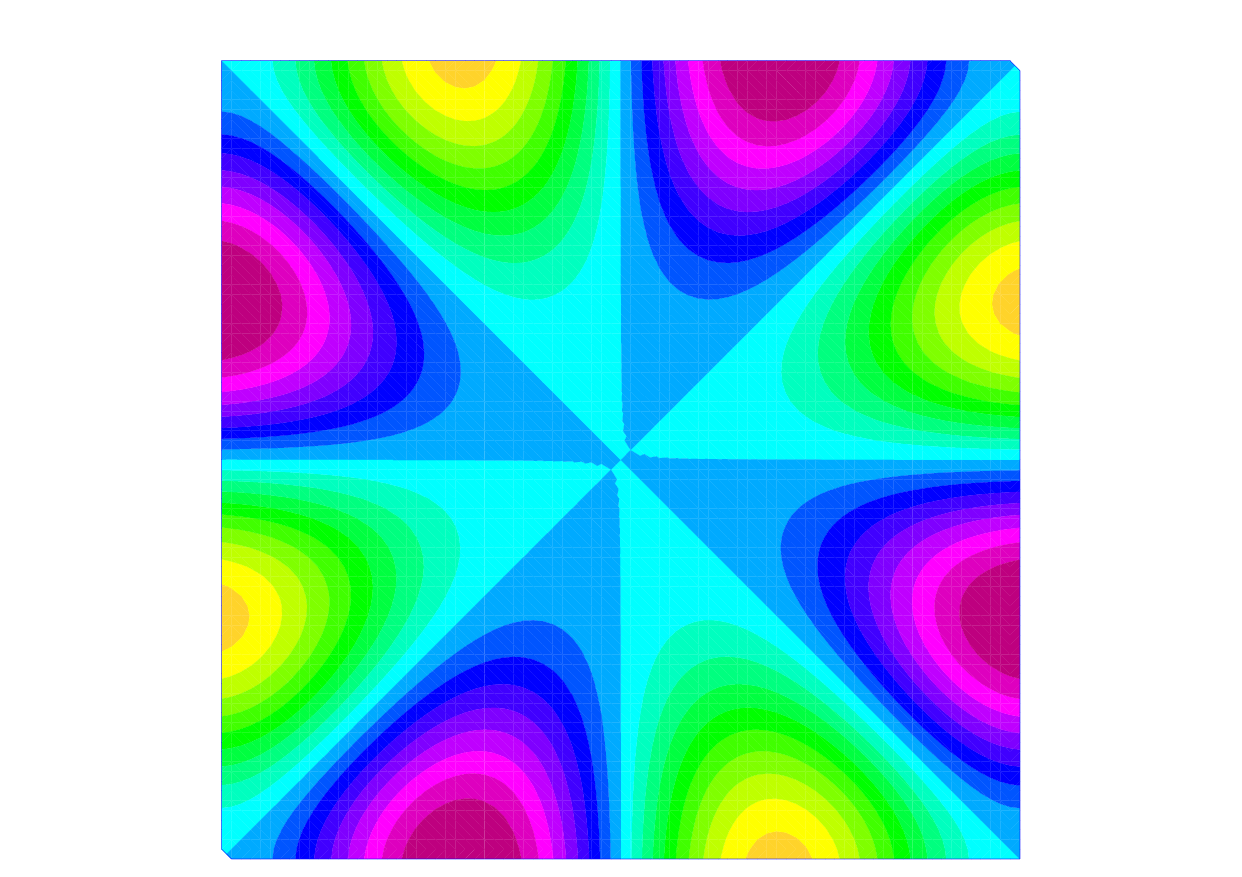}
	\caption{$\lambda_{9} = 4.6$}
	\label{SubFig:DtN3}
\end{subfigure}
\hfill
\begin{subfigure}[b]{0.18\textwidth}
	\centering
	\includegraphics[width=\textwidth,trim=7cm 0cm 7cm 0cm ,clip]{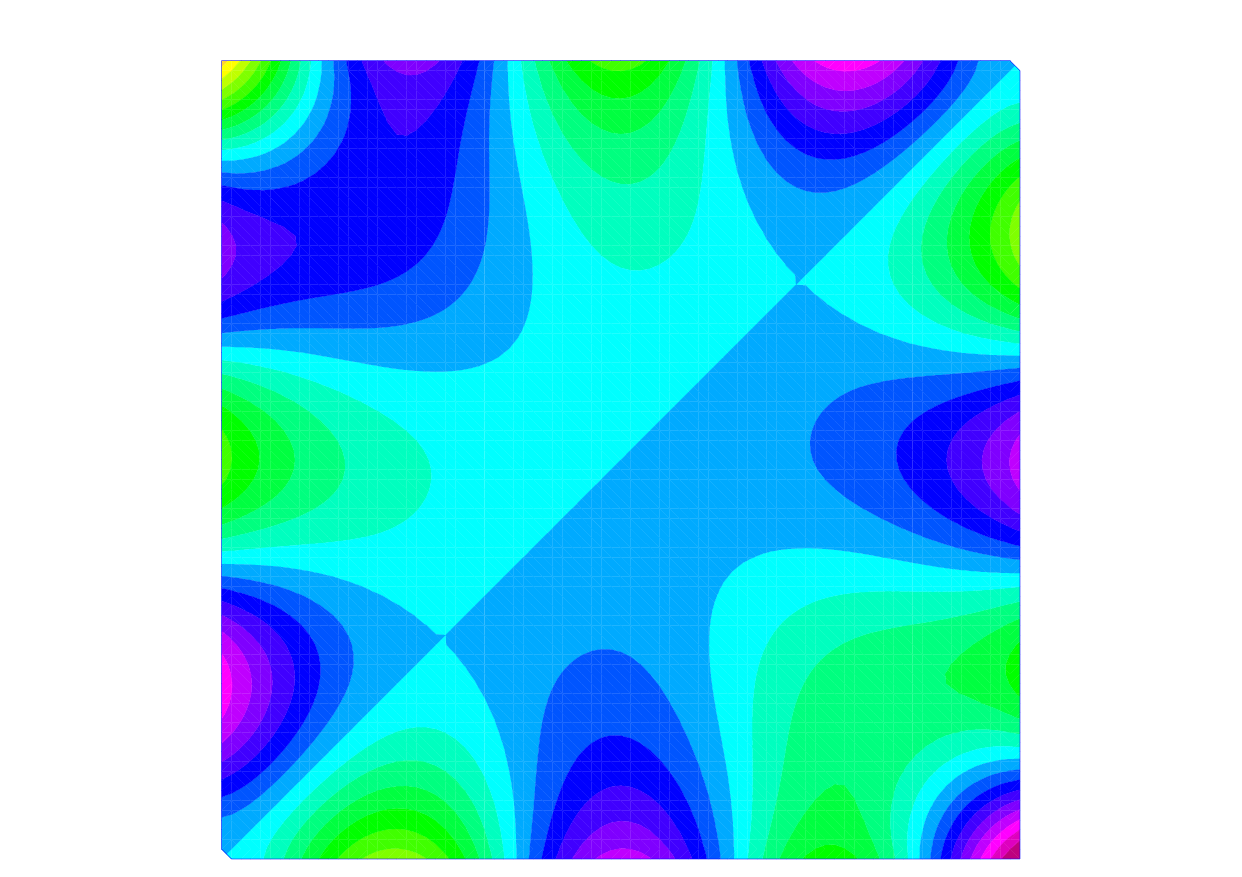}
	\caption{$\lambda_{13} = 32.0$}
\end{subfigure}
\hfill
\begin{subfigure}[b]{0.18\textwidth}
	\centering
	\includegraphics[width=\textwidth,trim=7cm 0cm 7cm 0cm ,clip]{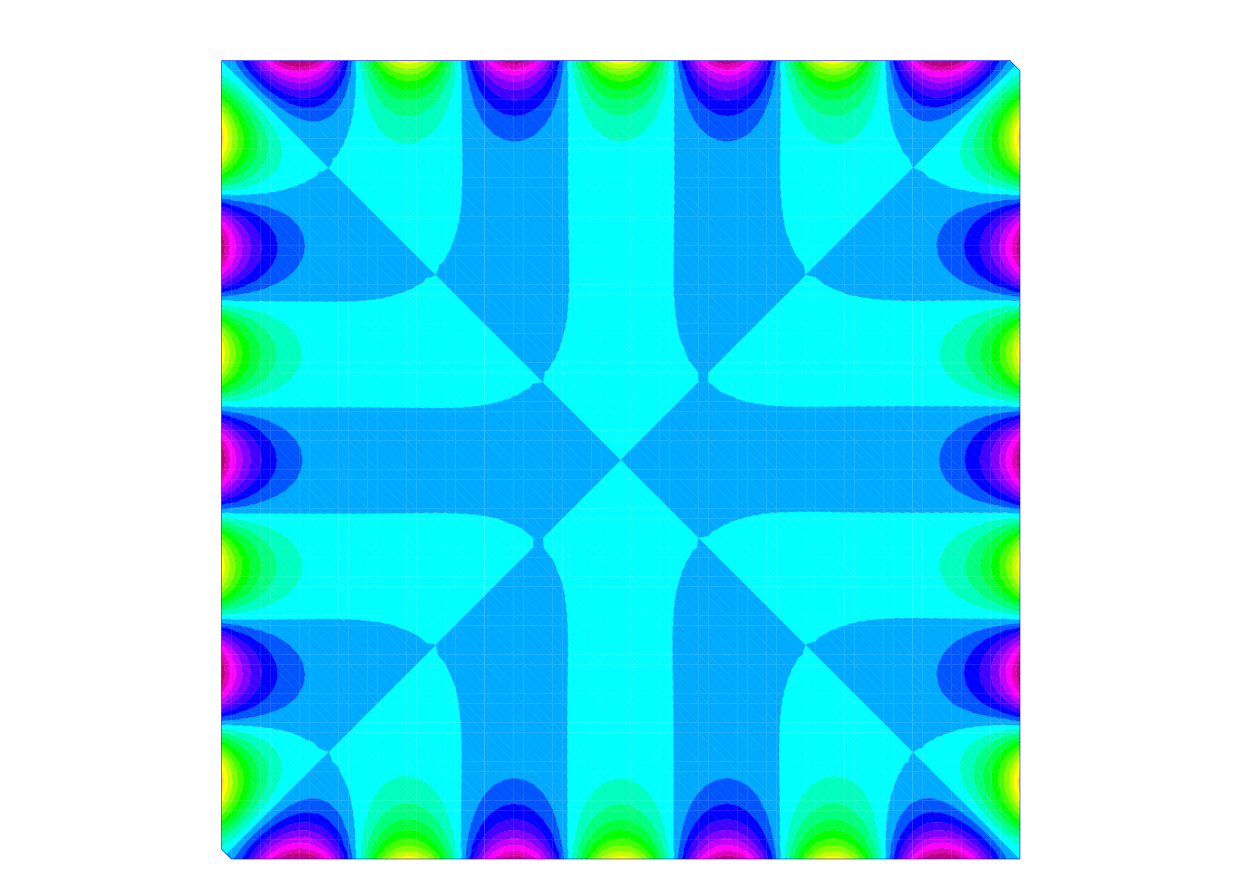}
	\caption{$\lambda_{28} = 107.5$}
	\label{SubFig:DtN5}
\end{subfigure}
\\[1ex]
\centering
\rotatebox{90}{\hspace*{1.2cm} H-GenEO}
\hfill
\begin{subfigure}[b]{0.18\textwidth}
	\centering
	\includegraphics[width=\textwidth,trim=7cm 0cm 7cm 0cm ,clip]{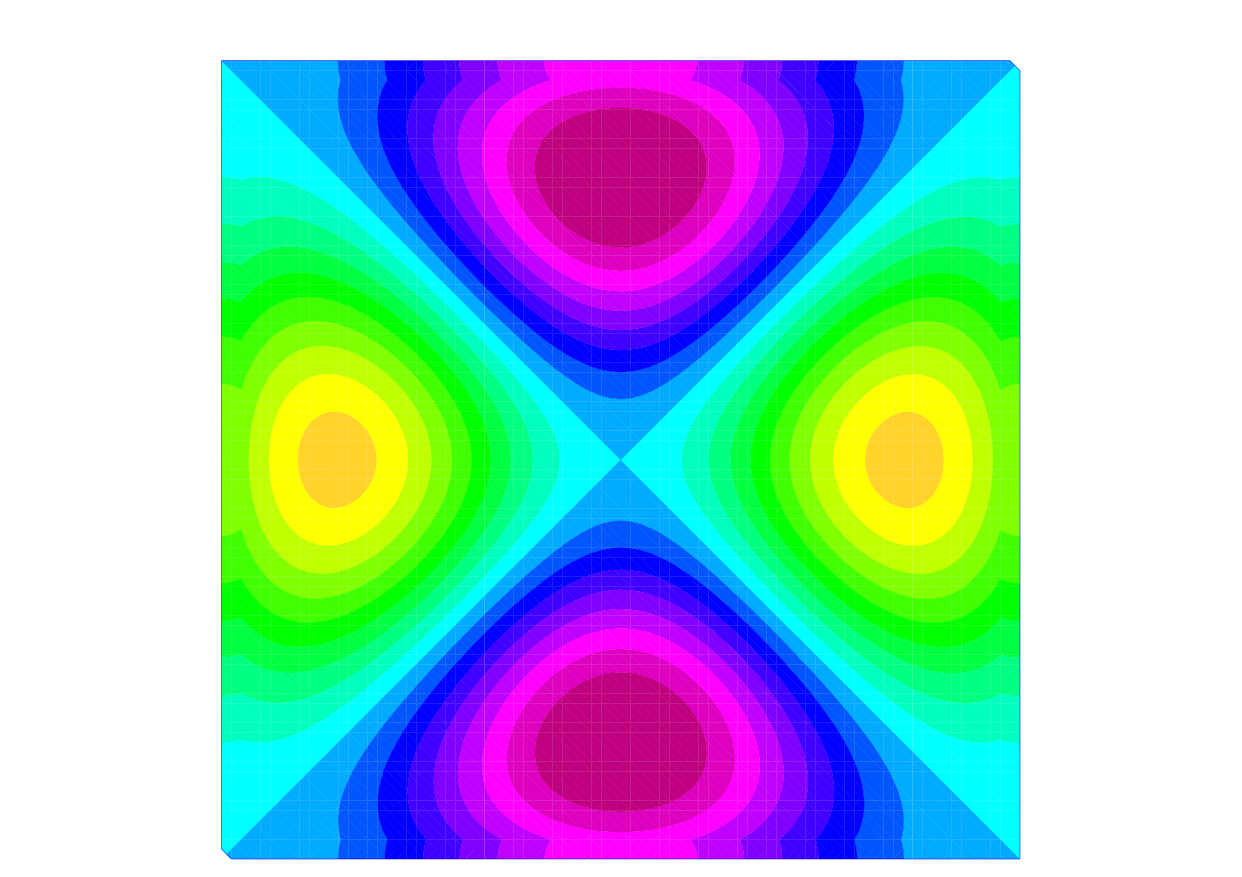}
	\caption{$\lambda_{3} = -0.258$}
\end{subfigure}
\hfill
\begin{subfigure}[b]{0.18\textwidth}
	\centering
	\includegraphics[width=\textwidth,trim=7cm 0cm 7cm 0cm ,clip]{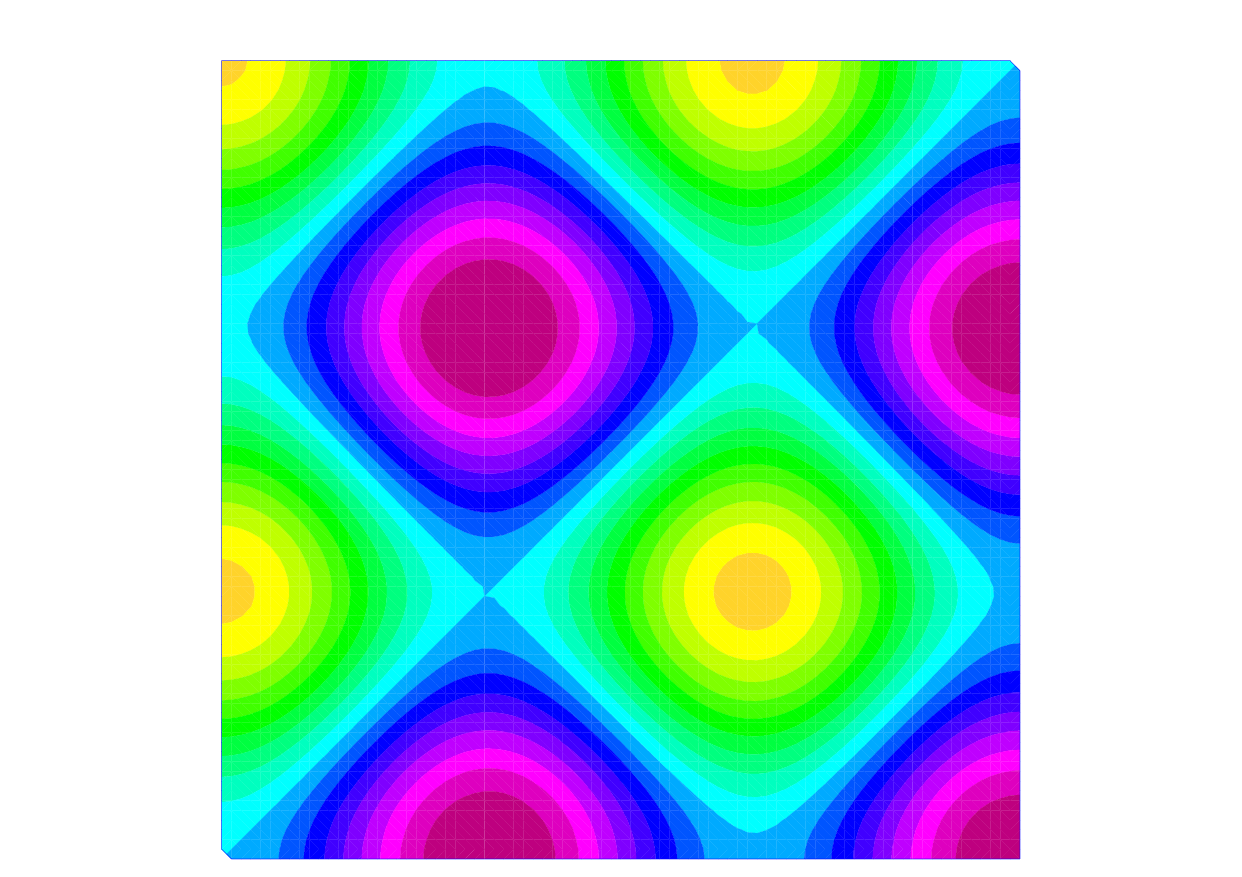}
	\caption{$\lambda_{7} = -0.006$}
\end{subfigure}
\hfill
\begin{subfigure}[b]{0.18\textwidth}
	\centering
	\includegraphics[width=\textwidth,trim=7cm 0cm 7cm 0cm ,clip]{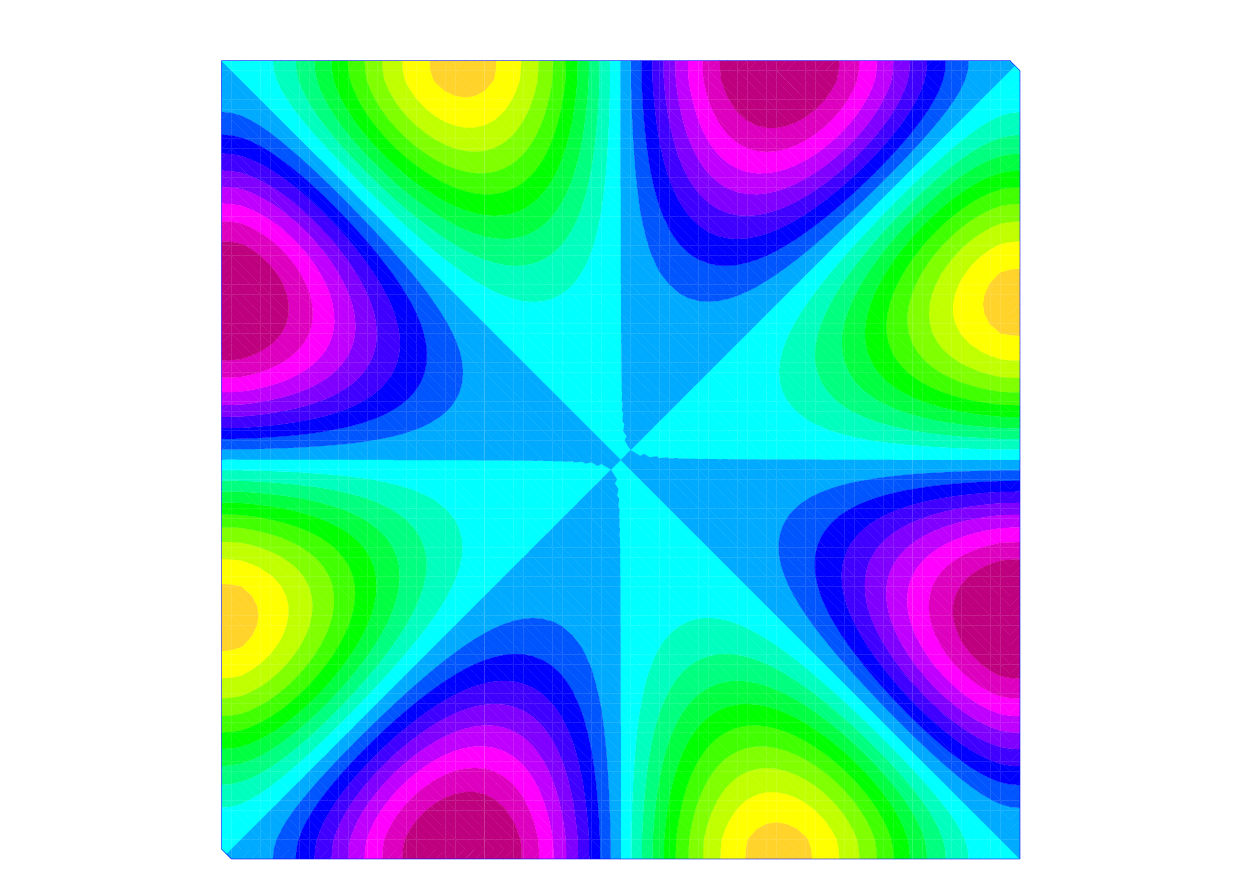}
	\caption{$\lambda_{9} = 0.019$}
\end{subfigure}
\hfill
\begin{subfigure}[b]{0.18\textwidth}
	\centering
	\includegraphics[width=\textwidth,trim=7cm 0cm 7cm 0cm ,clip]{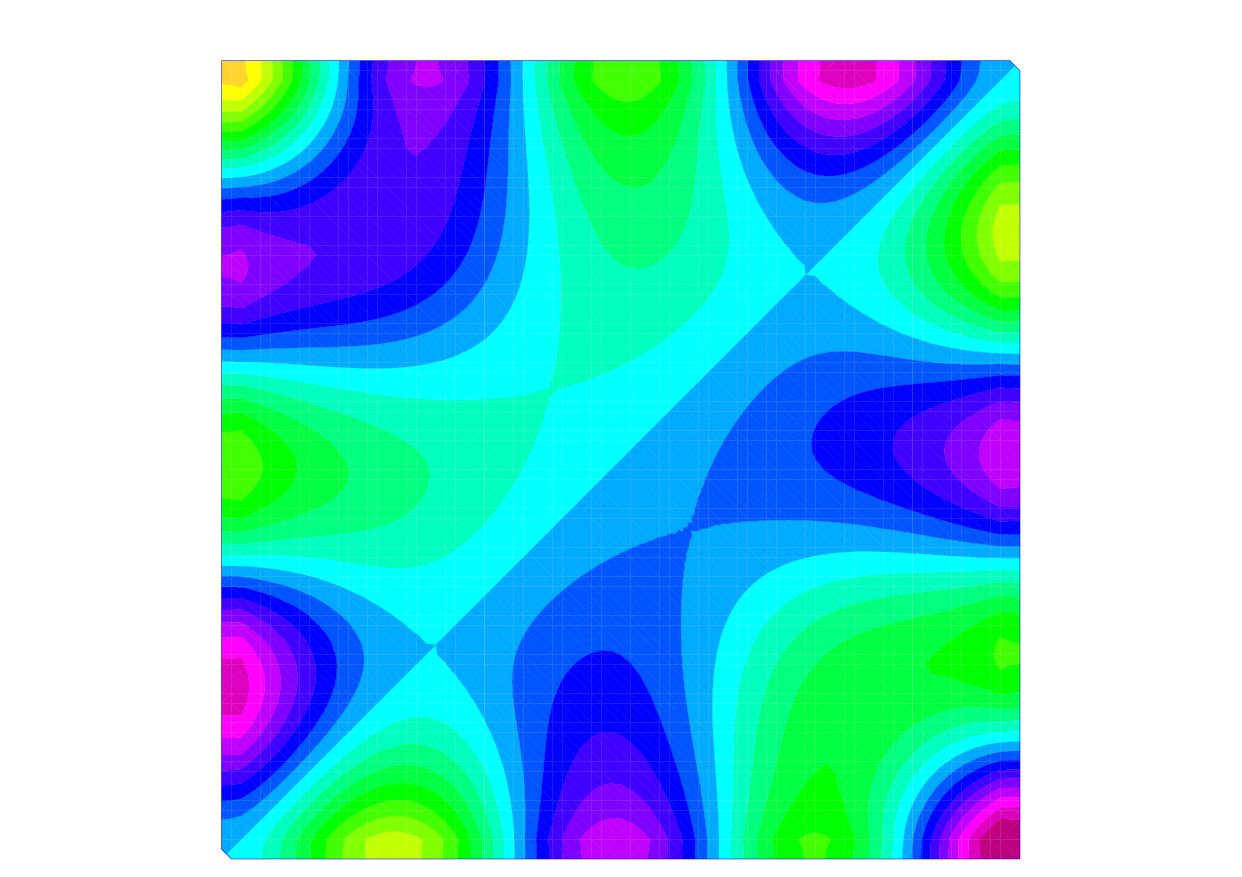}
	\caption{$\lambda_{14} = 0.149$}
\end{subfigure}
\hfill
\begin{subfigure}[b]{0.18\textwidth}
	\centering
	\includegraphics[width=\textwidth,trim=7cm 0cm 7cm 0cm ,clip]{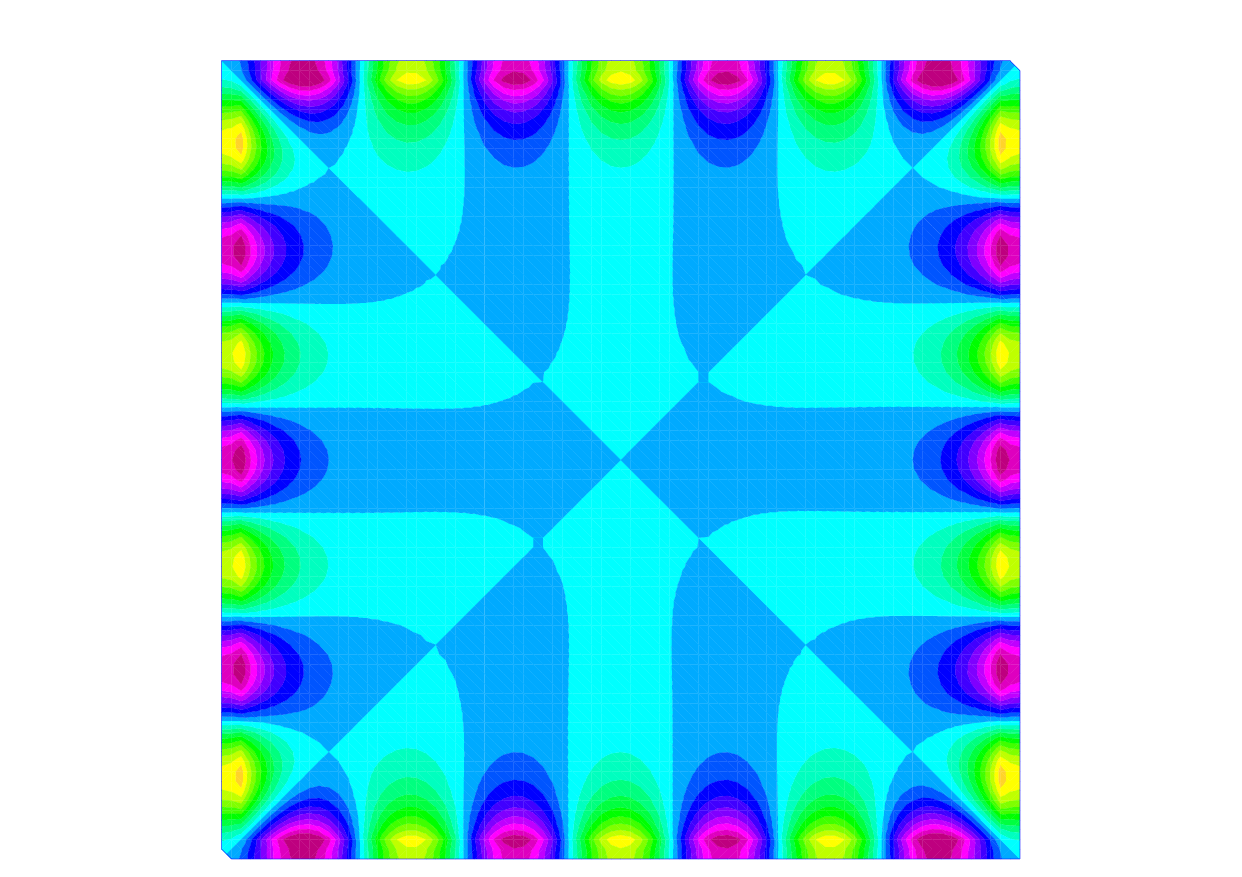}
	\caption{$\lambda_{34} = 0.466$}
\end{subfigure}
\\[1ex]
\centering
\rotatebox{90}{\hspace*{1.2cm} H-GenEO}
\hfill
\begin{subfigure}[b]{0.18\textwidth}
	\centering
	\includegraphics[width=\textwidth,trim=7cm 0cm 7cm 0cm ,clip]{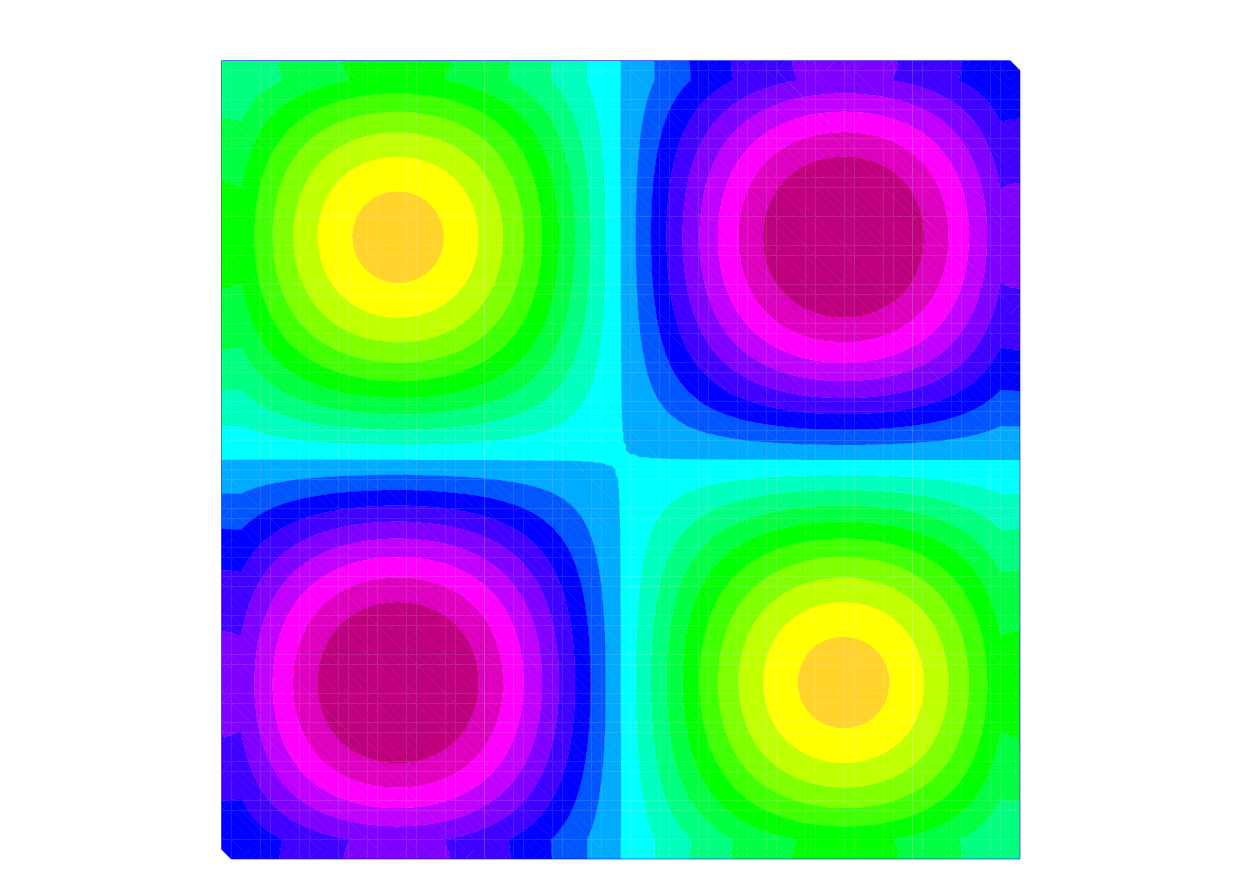}
	\caption{$\lambda_{1} = -0.431$}
\end{subfigure}
\hfill
\begin{subfigure}[b]{0.18\textwidth}
	\centering
	\includegraphics[width=\textwidth,trim=7cm 0cm 7cm 0cm ,clip]{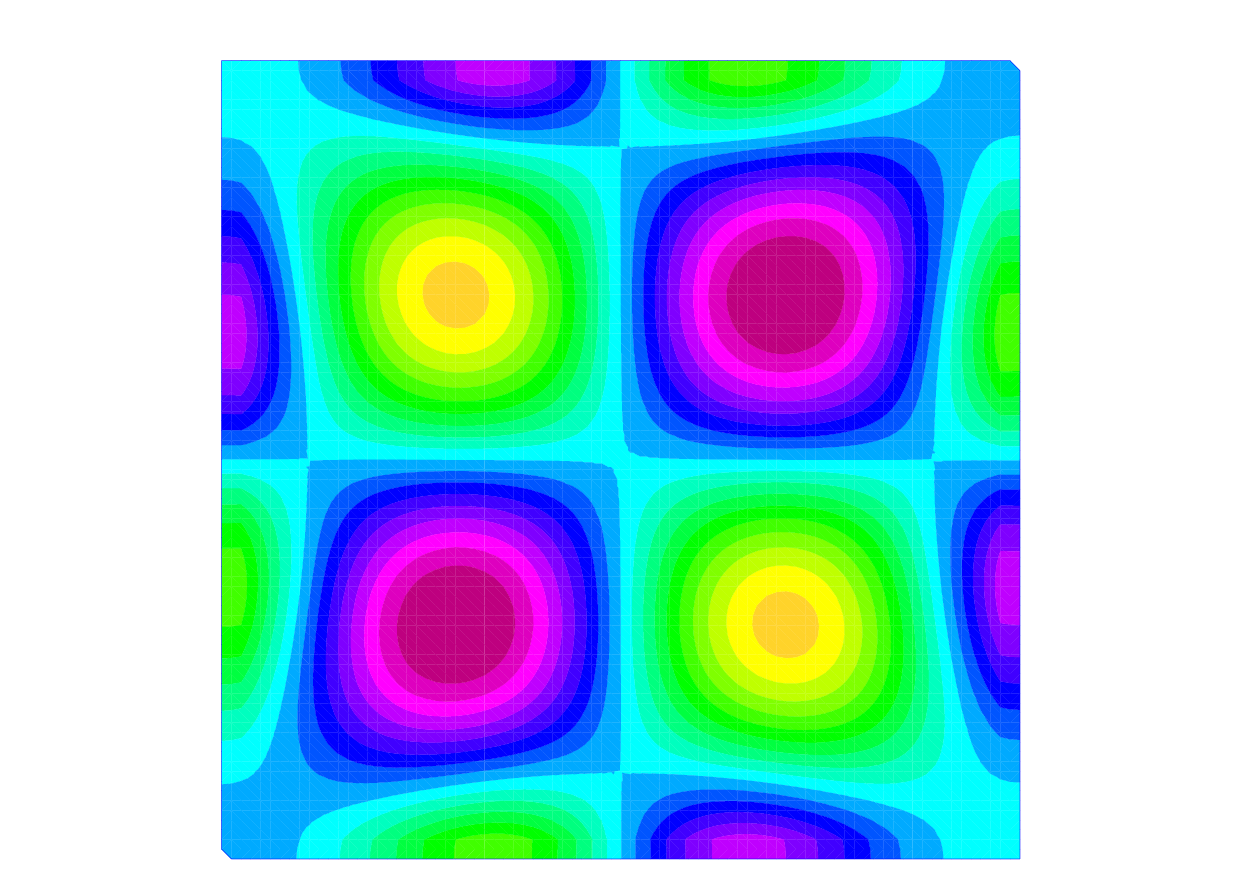}
	\caption{$\lambda_{17} = 0.233$}
\end{subfigure}
\hfill
\begin{subfigure}[b]{0.18\textwidth}
	\centering
	\includegraphics[width=\textwidth,trim=7cm 0cm 7cm 0cm ,clip]{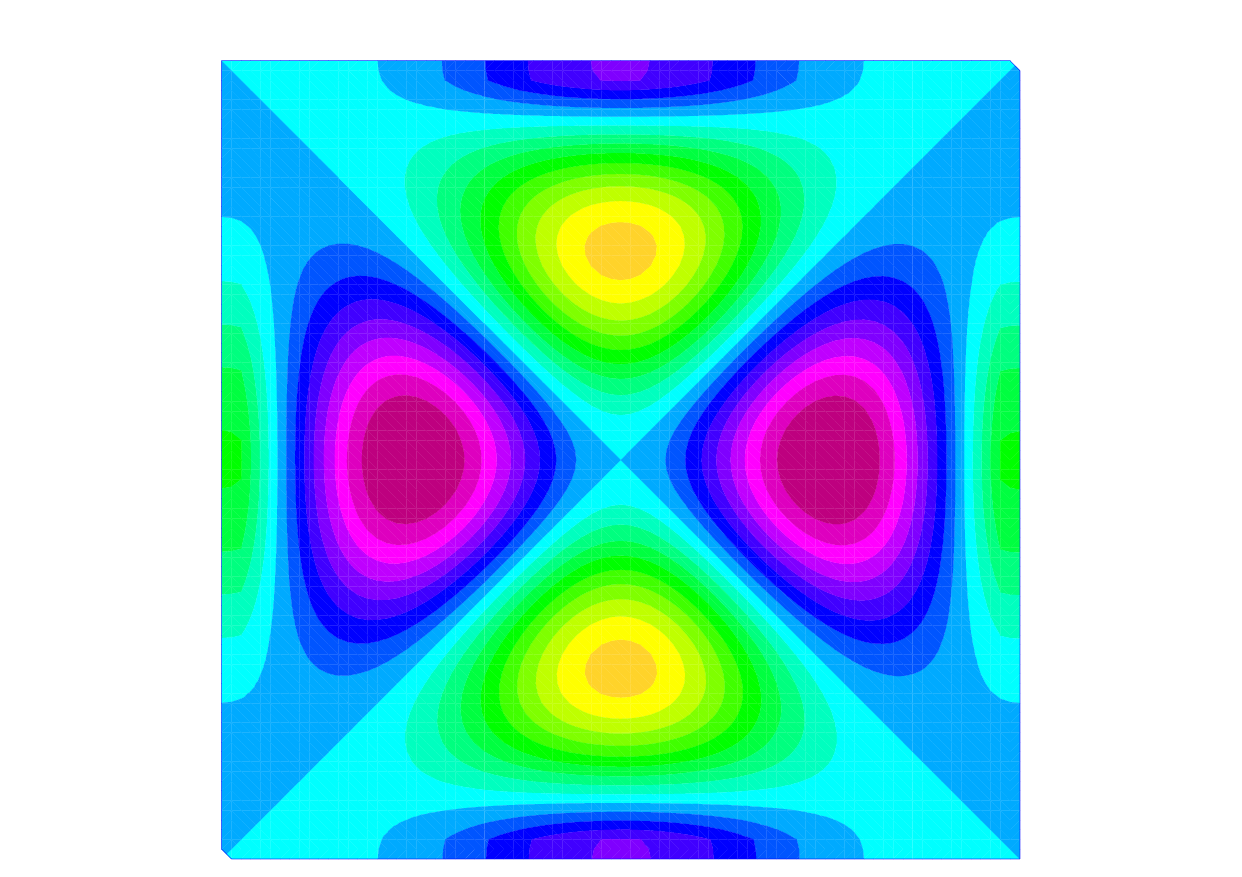}
	\caption{$\lambda_{22} = 0.315$}
\end{subfigure}
\hfill
\begin{subfigure}[b]{0.18\textwidth}
	\centering
	\includegraphics[width=\textwidth,trim=7cm 0cm 7cm 0cm ,clip]{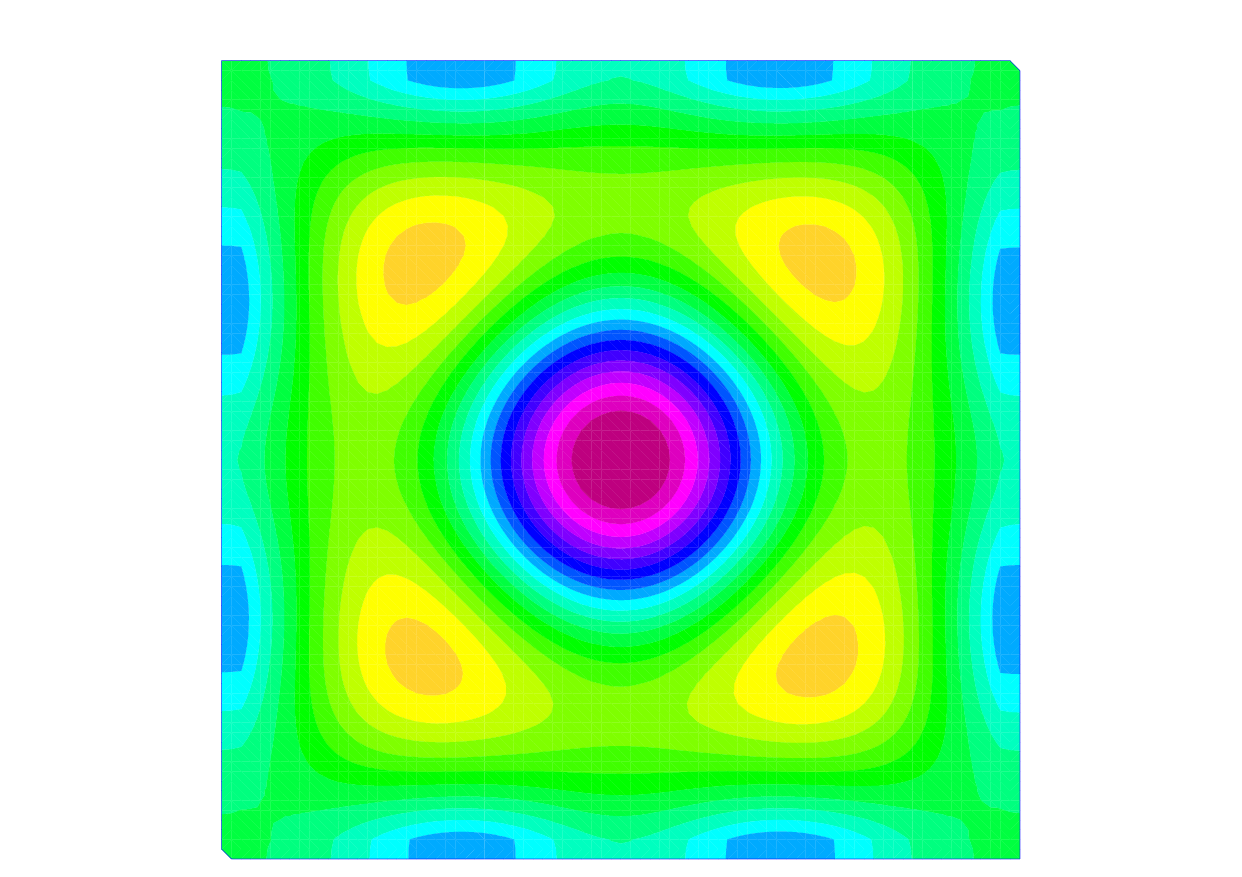}
	\caption{$\lambda_{23} = 0.315$}
\end{subfigure}
\hfill
\begin{subfigure}[b]{0.18\textwidth}
	\centering
	\includegraphics[width=\textwidth,trim=7cm 0cm 7cm 0cm ,clip]{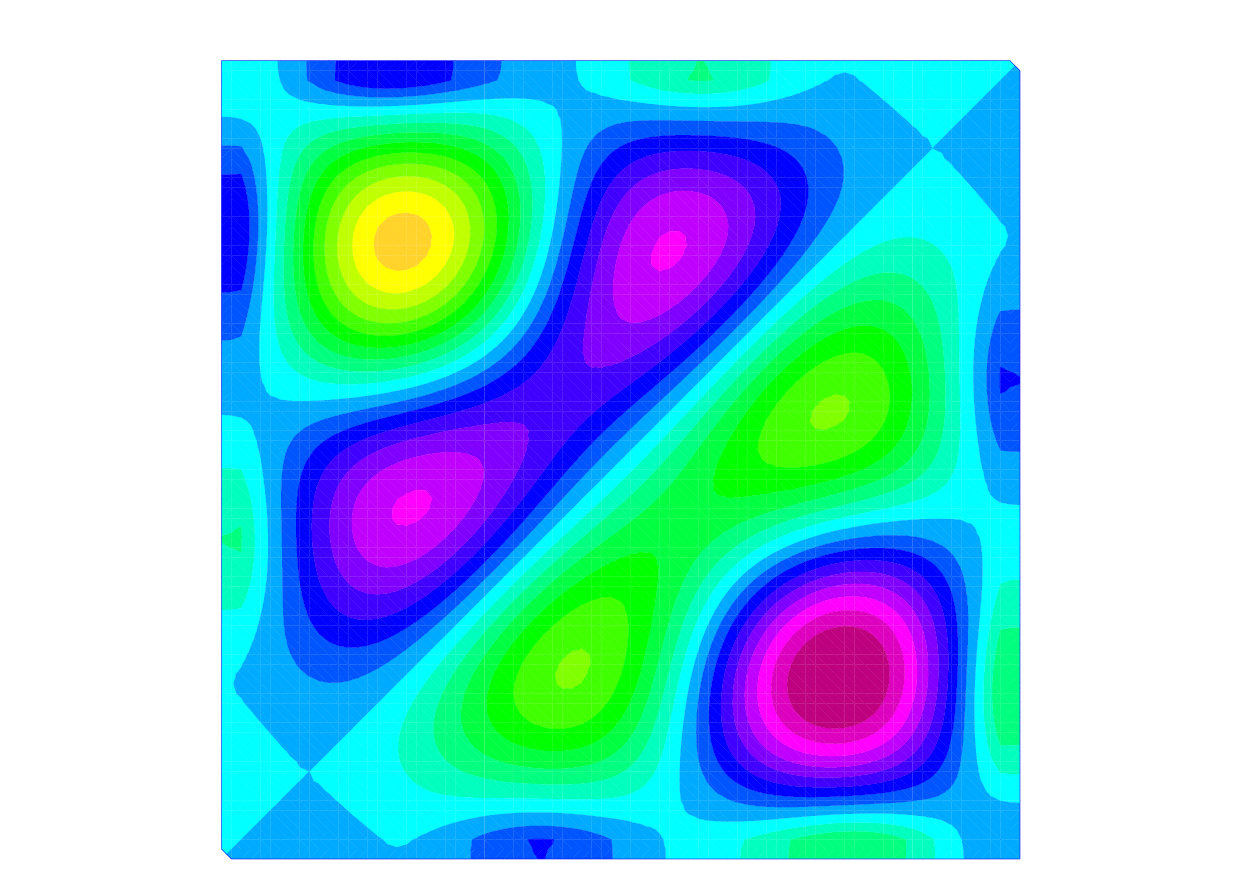}
	\caption{$\lambda_{33} = 0.438$}
\end{subfigure}
\caption{Local eigenfunctions for $k = 46.5$. \emph{Top row:} Examples using DtN \eqref{DiscreteDtNEigenproblem}. \emph{Middle row:} Equivalent examples using H-GenEO \eqref{DiscreteH-GenEOEigenproblem}. \emph{Bottom row:} Examples using H-GenEO which are not found amongst the DtN eigenfunctions.}
\label{Fig:Eigenfunctions}
\end{figure}

To further exhibit the link between these spectral coarse spaces for the Helmholtz problem, we consider example eigenfunctions for both the DtN \eqref{DiscreteDtNEigenproblem} and H-GenEO \eqref{DiscreteH-GenEOEigenproblem} eigenproblems on the central subdomain of a $5 \times 5$ decomposition into square subdomains for a homogeneous model problem with $k = 46.5$ (see Section~\ref{sec:Results} for details). Figure~\ref{Fig:Eigenfunctions} displays a selection of eigenfunctions and in the top two rows we plot DtN (top row) and H-GenEO (middle row) eigenfunctions which show the same features, visually being very similar. Note that, since the central subdomain does not touch the Robin boundary $\Gamma_{R}$, both eigenproblems are real and symmetric, albeit indefinite, and so all eigenvalues $\lambda$ are real. With DtN, as $\lambda$ increases, variation in the eigenfunctions tends to be restricted to the boundary, as can be seen in Figures~\ref{SubFig:DtN3}--\ref{SubFig:DtN5}. Such behaviour is observed for many H-GenEO eigenfunction too, however, we also obtain distinct eigenfunctions which are not found amongst DtN eigenfunctions: examples are given in the bottom row of Figure~\ref{Fig:Eigenfunctions} and we note that they tend to exhibit large variation in the interior of the subdomain. Nonetheless, there is a clear link between many of the DtN and H-GenEO eigenfunctions.

\begin{figure}[H]
	\widefigure
	\centering
	\hfill
	\begin{subfigure}[b]{0.49\textwidth}
		\centering
		\includegraphics[width=\textwidth,trim=0cm 0cm 0cm 0cm ,clip]{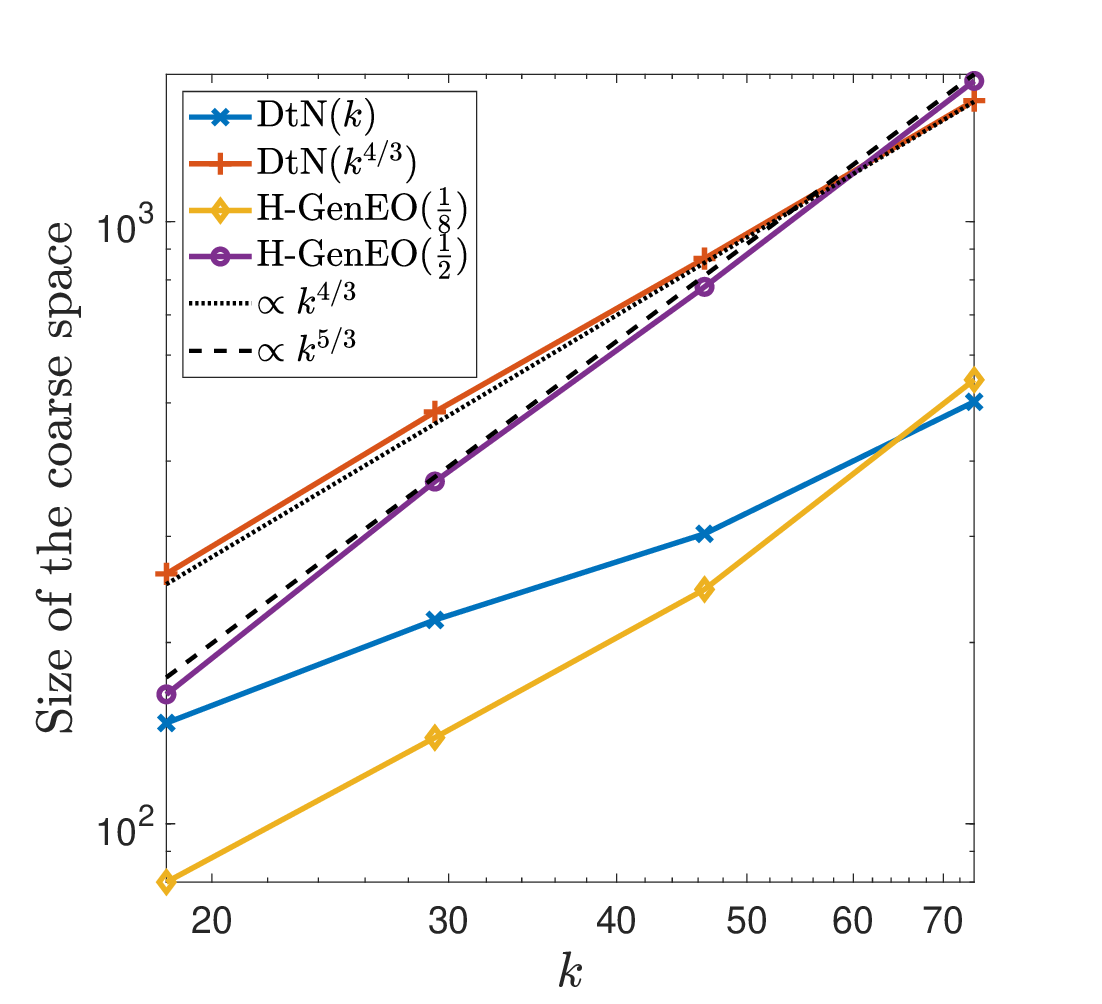}
		\caption{Varying the wave number $k$ for $N = 25$\\~\\}
		\label{Fig:CoarseSpaceSizesk}
	\end{subfigure}
	\hfill
	\begin{subfigure}[b]{0.49\textwidth}
		\centering
		\includegraphics[width=\textwidth,trim=0cm 0cm 0cm 0cm ,clip]{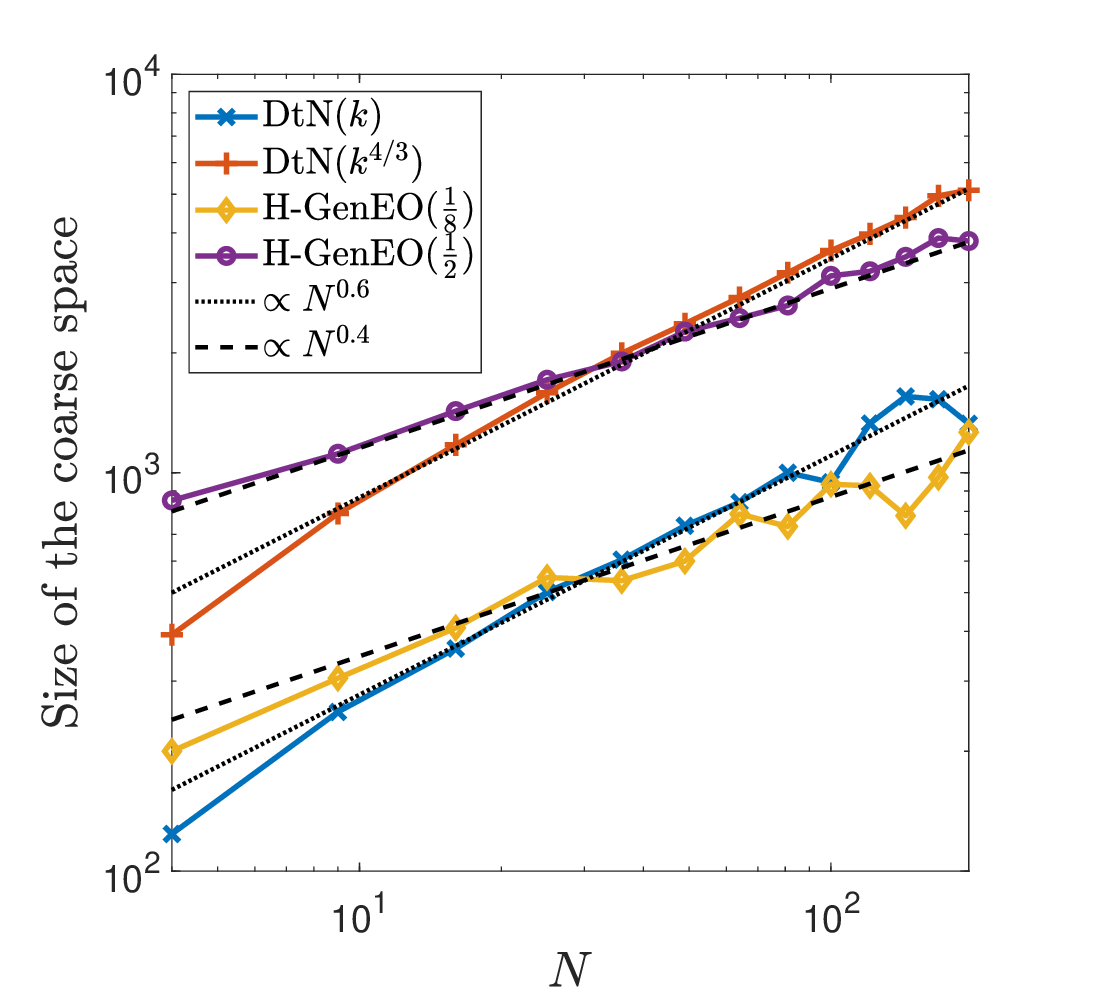}
		\caption{Varying the number of subdomains $N$ for $k = 73.8$\\~\\}
		\label{Fig:CoarseSpaceSizesN}
	\end{subfigure}
	\hfill
	\caption{The size of coarse space utilised for the homogeneous problem when using ORAS with the DtN and H-GenEO coarse spaces. A uniform decomposition into $\sqrt{N}\times\sqrt{N}$ square subdomains is used.}
	\label{Fig:CoarseSpaceSizes}
\end{figure}

\end{paracol}
\nointerlineskip
\begin{paracol}{2}
\switchcolumn

%%%%%%%%%%%%%%%%%%%%%%%%%%%%%%%%%%%%%%%%%%
\section{Results and Discussion}
\label{sec:Results}

In this section we present and discuss numerical results computed using FreeFEM \cite{Hecht:2012:NDI}, in particular through the functionality of \texttt{ffddm}, which handles the underlying domain decomposition data structures. As a model problem we consider the case of a wave guide in 2D, defined on the unit square $\Omega = (0,1)^2$. We impose homogeneous Dirichlet conditions on two opposite sides, namely \eqref{HelmholtzSpecificDirichletBC} with $u_{\Gamma_{D}} = 0$ on $\Gamma_{D} = \lbrace0,1\rbrace\times[0,1]$, and Robin conditions on the two remaining sides, that is \eqref{HelmholtzSpecificRobinBC} on $\Gamma_{R} = [0,1]\times\lbrace0,1\rbrace$. A point source is located in the centre of the domain at $(\frac{1}{2},\frac{1}{2})$ and provides the forcing function $f$. A schematic of this model problem is found in Figure~\ref{Fig:WaveGuide2D}.

\begin{figure}[H]
	\centering
	\begin{tikzpicture}[scale=1.5]
		% Grid
		\draw[step=0.5cm,gray!50,very thin, shift={(-1,-1)}] (0,0) grid (2,2);
		\draw[gray!50,very thin] (-1,-1) -- (1,1);
		\draw[gray!50,very thin] (-1,1) -- (1,-1);
		\draw[gray!50,very thin] (-1,0) -- (0,1) -- (1,0) -- (0,-1) -- cycle;
		% Boundaries
		\fill[gray,opacity=0.2] (-1,1) -- (-1,-1) -- (1,-1) -- (1,1) -- cycle;
		\draw (-1,1) -- (-1,-1) -- (1,-1) -- (1,1) -- cycle;
		% Boundary labels
		\draw (1.05,0) node[above,rotate=-90] {Dirichlet};
		\draw (-1.05,0) node[above,rotate=90] {Dirichlet};
		\draw (0,-1.05) node[below] {Robin};
		\draw (0,1.05) node[above] {Robin};
		% Coordinate labels
		\draw (-1.05,1) node[left] {$y = 1$}
		(-1.05,-1) node[left] {$y = 0$}
		(1,-1.05) node[below] {$x = 1$}
		(-1,-1.05) node[below] {$x = 0$};
		% Point source
		\draw[fill=black] (0,0) circle (0.025) node[below] {source $f$};
	\end{tikzpicture}
	\caption{\centering Schematic of the 2D wave guide model problem with example triangular mesh.}
	\label{Fig:WaveGuide2D}
\end{figure}
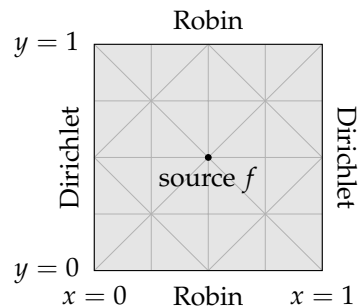

To discretise the problem we triangulate $\Omega$ using a Cartesian grid with spacing $h$ and alternating diagonals to form a simplicial mesh (see Figure~\ref{Fig:WaveGuide2D} where $h=\frac{1}{4}$). The discrete problem \eqref{LinearSystem} is then built from using P1 finite element approximation on this mesh. In order to avoid the pollution effect we choose $k$ and $h$ simultaneously so that $k^{3}h^{2} = \frac{2\pi}{10}$ is fixed. The large sparse linear system \eqref{LinearSystem} is solved using right-preconditioned GMRES, with the preconditioner given by the two-level ORAS method \eqref{2LevelAdaptiveDeflationPreconditioner} and choice of coarse space as stated. We terminate the GMRES iteration once a relative residual tolerance of $10^{-6}$ is reached. Unless otherwise stated, within the domain decomposition preconditioner we use minimal overlap: that is, one layer of adjoining mesh elements are added to the non-overlapping subdomains via the extension in \eqref{OverlappingSubdomains}. To solve the local eigenvalue problems in the two-level methods we make use of ARPACK \cite{Lehoucq:1998:ARPACK}, while both the subdomain solves and coarse space operator solves are given by MUMPS \cite{Amestoy:2001:AFA}.

We will first compare coarse spaces in the simplest case of a homogeneous problem, investigating the choice of eigenvalue threshold used. We then show results on how the methods perform in a variety of settings, for instance with non-uniform subdomains, heterogeneity and additional overlap. Finally, both weak and strong scalability tests are performed for H-GenEO applied to high wave number problems with timings reported.

\subsection{A Comparison of Methods for the Homogeneous Problem with Uniform Partitioning}

In Table~\ref{Table:WideMethodsN25} we give some benchmark results for the simplest problem of a homogeneous wave guide using a uniform decomposition into 25 square subdomains. We see that the one-level ORAS method \eqref{ORAS} performs relatively poorly as the wave number $k$ increases. The standard DtN coarse space \eqref{DiscreteDtNEigenproblem} (with $\eta_{\text{max}} = k$) is able to reduce iteration counts with a relatively small coarse space size, however, there is still a clear increase in iterations as $k$ increases. The $\Delta$-GenEO method \eqref{DiscreteLaplaceGenEOEigenproblem}, with $\lambda_{\text{max}} = \frac{1}{2}$, performs poorly here, often doing worse than the one-level method despite a larger coarse space than the DtN approach; this may be because the impedance conditions from the wave guide problem are not included in the definition of the $\Delta$-GenEO coarse space and so the eigenfunctions are not appropriate here. Finally, the standard H-GenEO method \eqref{DiscreteH-GenEOEigenproblem} (with $\eta_{\text{max}} = \frac{1}{2}$) performs well and significantly reduces the iteration counts, by a factor of 10 for the largest wave number, and provides robustness to increasing wave number $k$ (in fact iteration counts tend to decrease with $k$). We note that the size of the H-GenEO coarse space is larger than the DtN coarse space and we now explore this further.

\begin{specialtable}[H]
	\centering
	\caption{Preconditioned GMRES iteration counts and size of coarse space (in parentheses) for the homogeneous problem when using ORAS and various coarse spaces. A uniform decomposition into $5\times5$ square subdomains is used, giving 25 subdomains in total.}
	\label{Table:WideMethodsN25}
	\tabulinesep=1.2mm
	\begin{tabu}{cc|cccc}
		\toprule
		$k$ & $h^{-1}$ & one-level & DtN & $\Delta$-GenEO & H-GenEO \\
		\midrule
		18.5 & 100 &  73 & 19 (147) &  53  (135) & 21 (164) \\
		29.3 & 200 &  97 & 26 (218) & 100  (271) & 18 (370) \\
		46.5 & 400 & 125 & 35 (303) & 148  (560) & 17 (779) \\
		73.8 & 800 & 156 & 42 (502) & 220 (1120) & 15 (1712) \\
		\bottomrule
	\end{tabu}
\end{specialtable}

In Table~\ref{Table:DtNandH-GenEOThresholdsN25} we provide results for both the DtN and H-GenEO methods with differing eigenvalue thresholds $\eta_\text{max}$. For DtN we use the standard threshold $\eta_\text{max} = k$, the suggested threshold from \cite{Bootland:2019:ODN} $\eta_\text{max} = k^{4/3}$, and the larger threshold of $\eta_\text{max} = k^{3/2}$. For H-GenEO we use the standard threshold $\eta_\text{max} = \frac{1}{2}$ as well as the weaker thresholds of $\frac{1}{4}$ and $\frac{1}{8}$, the latter giving coarse space sizes more comparable to the standard DtN approach. To differentiate between these methods we use the notation DtN($\eta_\text{max}$) and H-GenEO($\eta_\text{max}$) where $\eta_\text{max}$ is as specified. We notice in Table~\ref{Table:DtNandH-GenEOThresholdsN25} that increasing the DtN threshold to $\eta_\text{max} = k^{4/3}$ significantly improves the iteration counts for this problem so that they are almost independent of $k$, albeit very slightly growing. Increasing the threshold further only marginally improves the iteration counts, which still grow slightly with $k$, but at the expense of a coarse space almost twice the size. On the other hand, if we relax the H-GenEO threshold we start to see higher iteration counts and lose some robustness but generally iteration counts do not increase with $k$ as they do for the standard DtN method. We note that, for this homogeneous problem, roughly comparable coarse space sizes give approximately similar iteration counts and so it is primarily the thresholds used in DtN and H-GenEO that dictate the different growth behaviour we observe.

\end{paracol}
\nointerlineskip
\begin{specialtable}[H]
	\widetable
	\centering
	\caption{Preconditioned GMRES iteration counts and size of coarse space (in parentheses) for the homogeneous problem when using ORAS and the DtN and H-GenEO coarse spaces with varying eigenvalue thresholds. A uniform decomposition into $5\times5$ square subdomains is used, giving 25 subdomains in total.}
	\label{Table:DtNandH-GenEOThresholdsN25}
	\tabulinesep=1.2mm
	\begin{tabu}{cc|ccc|ccc}
		\toprule
		$k$ & $h^{-1}$ & DtN($k$) & DtN($k^{4/3}$) & DtN($k^{3/2}$) & H-GenEO($\frac{1}{8}$) & H-GenEO($\frac{1}{4}$) & H-GenEO($\frac{1}{2}$) \\
		\midrule
		18.5 & 100 & 19 (147) & 13  (260) & 11  (403) & 46  (80) & 31 (105) & 21  (164) \\
		29.3 & 200 & 26 (218) & 14  (483) & 13  (759) & 53 (139) & 33 (189) & 18  (370) \\
		46.5 & 400 & 35 (303) & 14  (868) & 12 (1479) & 56 (245) & 35 (378) & 17  (779) \\
		73.8 & 800 & 42 (502) & 16 (1588) & 15 (2925) & 40 (546) & 25 (800) & 15 (1712) \\
		\bottomrule
	\end{tabu}
\end{specialtable}
\begin{paracol}{2}
\switchcolumn

To explore the growth in the size of the coarse space further, in Figure~\ref{Fig:CoarseSpaceSizesk} we plot coarse space size against the wave number $k$. From this we can see that growth for DtN($k^{4/3}$) is approximately proportional to $k^{4/3}$, while for H-GenEO($\frac{1}{2}$) it is around $k^{5/3}$ for our model problem. When the thresholds are relaxed it becomes less clear on the precise relationship but we note that the coarse space sizes grows more slowly with a weaker threshold, especially for the DtN approach. The faster growth seen for DtN($k^{4/3}$) and H-GenEO($\frac{1}{2}$) may help accommodate the stronger robustness to $k$ observed in the iteration counts of Table~\ref{Table:DtNandH-GenEOThresholdsN25}. These two approaches appear to provide the best trade-off for obtaining a well-behaved method and so we will focus primarily on these approaches, but first we consider the question of scalability.

\newpage
\end{paracol}
\nointerlineskip
\begin{paracol}{2}
\switchcolumn

\subsection{Scalability of DtN and H-GenEO for the Homogeneous Problem with Uniform Partitioning}

We now investigate the scalability of the DtN and H-GenEO methods. This will depend on the threshold used and so we compare results for DtN($k$) and H-GenEO($\frac{1}{8}$) as well as DtN($k^{4/3}$) and H-GenEO($\frac{1}{2}$), with each pair of approaches giving broadly similar iteration counts. Results for the homogeneous problem with $k = 73.8$ and $h^{-1} = 800$ are given in Table~\ref{Table:DtNandH-GenEO_Scalability} for an increasing number of uniform square subdomains $N$. We see that the DtN($k$) approach does not exhibit scalability here, with iteration counts that noticeably increase with $N$. Similarly, with the weaker threshold, H-GenEO($\frac{1}{8}$) also fails to be scalable. On the other hand, both DtN($k^{4/3}$) and H-GenEO($\frac{1}{2}$) are scalable here, with low iteration counts that vary little with $N$.

Comparing the size of the coarse spaces employed, we see that the DtN($k^{4/3}$) coarse space grows faster with $N$ and becomes larger than the H-GenEO($\frac{1}{2}$) coarse space, which may account for its particularly strong robustness to $N$ here. From Figure~\ref{Fig:CoarseSpaceSizesN} we see that for both DtN approaches the coarse space size grows approximately proportional to $N^{0.6}$ while for the H-GenEO method it is around $N^{0.4}$. This suggests H-GenEO may be advantageous when $N$ becomes large due to the smaller coarse space required. For each method, the average number of eigenvectors taken per subdomain decreases as $N$ increases, this means that as well as smaller subdomain solves we also benefit from requiring fewer eigenvectors to be computed per subdomain, even if the global coarse space size increases. We note that the size of the DtN($k$) or H-GenEO($\frac{1}{8}$) coarse space sometimes shrinks as we increase $N$ and in these cases the iteration counts tend to be particularly poor, suggesting that the thresholds of $\eta_\text{max} = k$ for DtN and $\eta_\text{max} = \frac{1}{8}$ are not doing a suitable job in capturing the eigenfunctions required for scalability. As such, we now narrow our focus to the DtN($k^{4/3}$) and H-GenEO($\frac{1}{2}$) approaches.

\newpage
\end{paracol}
\nointerlineskip
\begin{specialtable}[H]
\widetable
\centering
\caption{Preconditioned GMRES iteration counts (above), size of coarse space (middle), and average number of eigenvectors taken per subdomain (below) for the homogeneous problem when using ORAS with the DtN and H-GenEO coarse spaces and a varying number of subdomains $N$ for $k = 73.8$ and $h^{-1} = 800$. A uniform decomposition into $\sqrt{N}\times\sqrt{N}$ square subdomains is used.}
\label{Table:DtNandH-GenEO_Scalability}
\tabulinesep=1.2mm
\begin{tabu}{c|ccccccccccccc}
	\toprule
	$N$ & 4 & 9 & 16 & 25 & 36 & 49 & 64 & 81 & 100 & 121 & 144 & 169 & 196 \\
	\midrule
	DtN($k$) & 28 & 32 & 40 & 42 & 51 & 76 & 49 & 94 & 90 & 36 & 37 & 96 & 154 \\
	DtN($k^{4/3}$) & 15 & 16 & 19 & 16 & 16 & 16 & 15 & 16 & 15 & 15 & 16 & 17 & 17 \\
	H-GenEO($\frac{1}{8}$) & 26 & 31 & 36 & 40 & 71 & 70 & 65 & 127 & 81 & 116 & 247 & 194 & 138 \\
	H-GenEO($\frac{1}{2}$) & 13 & 15 & 15 & 15 & 16 & 16 & 16 & 18 & 16 & 18 & 18 & 18 & 19 \\
	\midrule\midrule
	DtN($k$) & 124 & 251 & 362 & 502 & 605 & 736 & 843 & 1000 & 946 & 1329 & 1554 & 1529 & 1327 \\
	DtN($k^{4/3}$) & 392 & 790 & 1175 & 1588 & 1994 & 2366 & 2753 & 3176 & 3611 & 3976 & 4369 & 4955 & 5118 \\
	H-GenEO($\frac{1}{8}$) & 200 & 305 & 408 & 546 & 536 & 600 & 788 & 733 & 936 & 927 & 780 & 974 & 1264 \\
	H-GenEO($\frac{1}{2}$) & 852 & 1116 & 1428 & 1712 & 1903 & 2261 & 2444 & 2629 & 3120 & 3204 & 3482 & 3882 & 3816 \\
	\midrule\midrule
	DtN($k$) & 31.0 & 27.9 & 22.6 & 20.1 & 16.8 & 15.0 & 13.2 & 12.3 &  9.5 & 11.0 & 10.8 &  9.0 &  6.8 \\
	DtN($k^{4/3}$) & 98.0 & 87.8 & 73.4 & 63.5 & 55.4 & 48.3 & 43.0 & 39.2 & 36.1 & 32.9 & 30.3 & 29.3 & 26.1 \\
	H-GenEO($\frac{1}{8}$) & 50.0 & 33.9 & 25.5 & 21.8 & 14.9 & 12.2 & 12.3 &  9.0 &  9.4 &  7.7 &  5.4 &  5.8 &  6.4 \\
	H-GenEO($\frac{1}{2}$) & 213.0 & 124.0 & 89.3 & 68.5 & 52.9 & 46.1 & 38.2 & 32.5 & 31.2 & 26.5 & 24.2 & 23.0 & 19.5 \\
	\bottomrule
\end{tabu}
\end{specialtable}
\begin{paracol}{2}
\switchcolumn

\subsection{Robustness of DtN and H-GenEO for the Homogeneous Problem with METIS Decomposition}

We now consider utilising non-uniform subdomains, as provided through the software METIS \cite{Karypis:1998:AFA}. In Table~\ref{Table:DtNandH-GenEO_METIS} we compare DtN($k^{4/3}$) and H-GenEO($\frac{1}{2}$) in this situation, again for the homogeneous problem. We observe that both methods retain their robustness, meaning that iteration counts only depend mildly on the wave number $k$ and the number of subdomains $N$. For H-GenEO($\frac{1}{2}$) the scalability becomes more favourable for larger $k$ and we can see that for the smallest wave numbers with large $N$ the average number of eigenvectors per subdomain becomes very small, in fact on many subdomains only a single eigenvector is taken and so the achieved tolerance on the eigenvalue may be somewhat weaker than $\frac{1}{2}$, which may explain the slightly poorer performance. One way this could be overcome is by always taking a minimum number of eigenvectors per subdomain; for instance, using at least 5 eigenvectors per subdomain when $k = 18.5$ gives iteration counts bounded by 19. Since we are primarily interested in approaches that remain effective for increasingly large wave numbers, where this issue diminishes, we do not worry further about this and assert that for problems of interest H-GenEO($\frac{1}{2}$) provides good scalability. On the other hand, for DtN($k^{4/3}$) the mild increase in iteration count is seen as $k$ increases here, with scalability observed for all wave numbers.

We note that the coarse space sizes for each method tends to be slightly larger with the more general decompositions used by METIS but otherwise the same trends are seen. As such, we conclude that non-uniform decompositions can be well-handled by the spectral coarse spaces employed here.

\newpage
\end{paracol}
\nointerlineskip
\begin{specialtable}[H]
	\widetable
	\centering
	\caption{Preconditioned GMRES iteration counts (above), size of coarse space (middle), and average number of eigenvectors taken per subdomain (below) for the homogeneous problem when using ORAS with DtN($k^{4/3}$) or H-GenEO($\frac{1}{2}$) and a varying number of subdomains. A non-uniform decomposition into $N$ subdomains is used, given by METIS.}
	\label{Table:DtNandH-GenEO_METIS}
	\tabulinesep=1.2mm
	\begin{tabu}{cc|cccccc|cccccc}
		\toprule
		& & \multicolumn{12}{c}{Number of subdomains $N$} \\
		& & \multicolumn{6}{c}{DtN($k^{4/3}$)} & \multicolumn{6}{c}{H-GenEO($\frac{1}{2}$)} \\
		$k$ & $h^{-1}$ & 20 & 40 & 80 & 120 & 160 & 200 & 20 & 40 & 80 & 120 & 160 & 200 \\
		\midrule
		18.5 & 100 & 10 & 10 & 10 & 10 & 10 & 10 & 15 & 17 & 19 & 22 & 27 & 27 \\
		29.3 & 200 & 12 & 15 & 11 & 12 & 12 & 12 & 15 & 17 & 19 & 20 & 22 & 23 \\
		46.5 & 400 & 12 & 13 & 15 & 13 & 13 & 13 & 15 & 16 & 16 & 18 & 20 & 20 \\
		73.8 & 800 & 15 & 15 & 14 & 16 & 14 & 16 & 15 & 16 & 17 & 17 & 17 & 19 \\
		117.2 & 1600 & 14 & 15 & 16 & 17 & 15 & 16 & 14 & 15 & 15 & 16 & 16 & 16 \\
		\midrule\midrule
		18.5 & 100 & 281 & 422 & 652 & 843 & 1005 & 1157 & 201 & 285 & 383 & 471 & 524 & 589 \\
		29.3 & 200 & 477 & 758 & 1130 & 1410 & 1693 & 1922 & 400 & 574 & 783 & 958 & 1097 & 1245 \\
		46.5 & 400 & 959 & 1466 & 2132 & 2677 & 3151 & 3553 & 869 & 1193 & 1670 & 2008 & 2253 & 2507 \\
		73.8 & 800 & 1695 & 2563 & 3751 & 4672 & 5486 & 6199 & 1863 & 2456 & 3433 & 4147 & 4749 & 5338 \\
		117.2 & 1600 & 3049 & 4695 & 6831 & 8486 & 9896 & 11092 & 4238 & 5680 & 7575 & 9049 & 10273 & 11305 \\
		\midrule\midrule
		18.5 & 100 & 14.1 & 10.6 &  8.2 &  7.0 &  6.3 &  5.8 & 10.1 &  7.1 &  4.8 &  3.9 &  3.3 &  2.9 \\
		29.3 & 200 & 23.9 & 18.9 & 14.1 & 11.8 & 10.6 &  9.6 & 20.0 & 14.3 &  9.8 &  8.0 &  6.9 &  6.2 \\
		46.5 & 400 & 48.0 & 36.6 & 26.6 & 22.3 & 19.7 & 17.8 & 43.5 & 29.8 & 20.9 & 16.7 & 14.1 & 12.5 \\
		73.8 & 800 & 84.8 & 64.1 & 46.9 & 38.9 & 34.3 & 31.0 & 93.2 & 61.4 & 42.9 & 34.6 & 29.7 & 26.7 \\
		117.2 & 1600 & 152.4 & 117.4 & 85.4 & 70.7 & 61.9 & 55.5 & 211.9 & 142.0 & 94.7 & 75.4 & 64.2 & 56.5 \\
		\bottomrule
	\end{tabu}
\end{specialtable}
\begin{paracol}{2}
\switchcolumn

\subsection{The Effect of Heterogeneity}

We now turn our attention to the key property of robustness to heterogeneities. For this we consider layered media within the wave guide. Three configurations, each having ten layers, will be used and are detailed in Figure~\ref{Fig:Layers}. The heterogeneity is introduced in the wave speed $c(\boldsymbol{x})$ and in each case $c$ takes values from $1$ to $\rho$, where $\rho$ is a contrast parameter determining the strength of the heterogeneity. The wave number is then given by $k = \omega/c$ where $\omega$ is the angular frequency; we will vary both $\omega$ and $\rho$ in our tests. Note that for the DtN method the eigenvalue threshold will now depend on $k_{s} = \max_{\vec{x}\in\Omega_{s}} k(\vec{x})$ which may be different for different subdomains $\Omega_{s}$. To avoid notational clutter we omit the subscript when referring to the method, namely retaining the name DtN($k^{4/3}$).

\begin{figure}[H]
	\centering
	\hspace*{\fill}
	\subfloat[\centering Increasing layers]{
		\begin{tikzpicture}[scale=1.5]
			\fill[gray,opacity=0.1] (-1,1) -- (1,1) -- (1,0.8) -- (-1,0.8) -- cycle;
			\fill[gray,opacity=0.2] (-1,0.6) -- (1,0.6) -- (1,0.8) -- (-1,0.8) -- cycle;
			\fill[gray,opacity=0.3] (-1,0.6) -- (1,0.6) -- (1,0.4) -- (-1,0.4) -- cycle;
			\fill[gray,opacity=0.4] (-1,0.2) -- (1,0.2) -- (1,0.4) -- (-1,0.4) -- cycle;
			\fill[gray,opacity=0.5] (-1,0.2) -- (1,0.2) -- (1,0.0) -- (-1,0.0) -- cycle;
			\fill[gray,opacity=0.6] (-1,-0.2) -- (1,-0.2) -- (1,0.0) -- (-1,0.0) -- cycle;
			\fill[gray,opacity=0.7] (-1,-0.2) -- (1,-0.2) -- (1,-0.4) -- (-1,-0.4) -- cycle;
			\fill[gray,opacity=0.8] (-1,-0.6) -- (1,-0.6) -- (1,-0.4) -- (-1,-0.4) -- cycle;
			\fill[gray,opacity=0.9] (-1,-0.6) -- (1,-0.6) -- (1,-0.8) -- (-1,-0.8) -- cycle;
			\fill[gray,opacity=1.0] (-1,-1) -- (1,-1) -- (1,-0.8) -- (-1,-0.8) -- cycle;
			\draw (-1,1) -- (-1,-1) -- (1,-1) -- (1,1) -- cycle;
			\label{Fig:IncreasingLayers}
		\end{tikzpicture}
	} \hspace*{\fill}
	\subfloat[\centering Alternating layers]{
		\begin{tikzpicture}[scale=1.5]
			\fill[gray,opacity=0.1] (-1,1) -- (1,1) -- (1,0.8) -- (-1,0.8) -- cycle;
			\fill[gray,opacity=1.0] (-1,0.6) -- (1,0.6) -- (1,0.8) -- (-1,0.8) -- cycle;
			\fill[gray,opacity=0.1] (-1,0.6) -- (1,0.6) -- (1,0.4) -- (-1,0.4) -- cycle;
			\fill[gray,opacity=1.0] (-1,0.2) -- (1,0.2) -- (1,0.4) -- (-1,0.4) -- cycle;
			\fill[gray,opacity=0.1] (-1,0.2) -- (1,0.2) -- (1,0.0) -- (-1,0.0) -- cycle;
			\fill[gray,opacity=1.0] (-1,-0.2) -- (1,-0.2) -- (1,0.0) -- (-1,0.0) -- cycle;
			\fill[gray,opacity=0.1] (-1,-0.2) -- (1,-0.2) -- (1,-0.4) -- (-1,-0.4) -- cycle;
			\fill[gray,opacity=1.0] (-1,-0.6) -- (1,-0.6) -- (1,-0.4) -- (-1,-0.4) -- cycle;
			\fill[gray,opacity=0.1] (-1,-0.6) -- (1,-0.6) -- (1,-0.8) -- (-1,-0.8) -- cycle;
			\fill[gray,opacity=1.0] (-1,-1) -- (1,-1) -- (1,-0.8) -- (-1,-0.8) -- cycle;
			\draw (-1,1) -- (-1,-1) -- (1,-1) -- (1,1) -- cycle;
			\label{Fig:AlternatingLayers}
		\end{tikzpicture}
	} \hspace*{\fill}
	\subfloat[\centering Diagonal layers]{
		\begin{tikzpicture}[scale=1.5]
			\fill[gray,opacity=1.0] (1,-1) -- (0.6,-1) -- (1,-0.6) -- cycle;
			\fill[gray,opacity=0.6] (0.6,-1) -- (0.2,-1) -- (1,-0.2) -- (1,-0.6) -- cycle;
			\fill[gray,opacity=1.0] (0.2,-1) -- (-0.2,-1) -- (1,0.2) -- (1,-0.2) -- cycle;
			\fill[gray,opacity=0.2] (-0.2,-1) -- (-0.6,-1) -- (1,0.6) -- (1,0.2) -- cycle;
			\fill[gray,opacity=1.0] (-0.6,-1) -- (-1,-1) -- (1,1) -- (1,0.6) -- cycle;
			\fill[gray,opacity=0.05] (-1,-1) -- (-1,-0.6) -- (0.6,1) -- (1,1) -- cycle;
			\fill[gray,opacity=1.0] (-1,-0.6) -- (-1,-0.2) -- (0.2,1) -- (0.6,1) -- cycle;
			\fill[gray,opacity=0.4] (-1,-0.2) -- (-1,0.2) -- (-0.2,1) -- (0.2,1) -- cycle;
			\fill[gray,opacity=1.0] (-1,0.2) -- (-1,0.6) -- (-0.6,1) -- (-0.2,1) -- cycle;
			\fill[gray,opacity=0.8] (-1,0.6) -- (-1,1) -- (-0.6,1) -- cycle;
			\draw (-1,1) -- (-1,-1) -- (1,-1) -- (1,1) -- cycle;
			\label{Fig:DiagonalLayers}
		\end{tikzpicture}
	}
	\hspace*{\fill}
	\caption{Piecewise constant layer profiles for the wave speed $c(\boldsymbol{x})$. For the darkest shade $c(\boldsymbol{x}) = 1$ while for the lightest shade $c(\boldsymbol{x}) = \rho$, with $\rho$ being the contrast factor.}
	\label{Fig:Layers}
\end{figure}
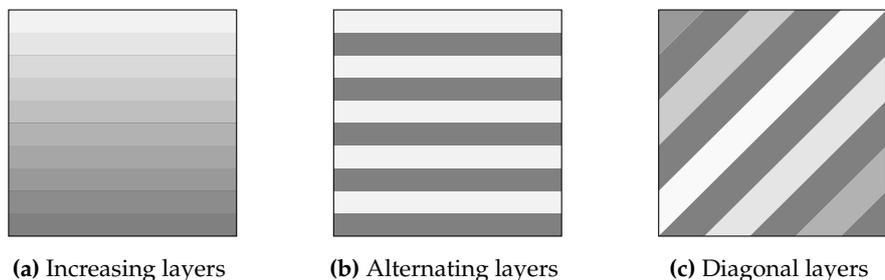

\newpage
Results for DtN($k^{4/3}$) and H-GenEO($\frac{1}{2}$) for the increasing layers problem (Figure~\ref{Fig:IncreasingLayers}) are provided in Table~\ref{Table:DtNandH-GenEO_IncreasingLayers}. Unfortunately we see that the DtN($k^{4/3}$) loses robustness with the heterogeneity present in this problem. In particular, we lose any robustness to the wave number $k$ and for the largest wave number used we also see that changes in the contrast, given by $\rho$, can begin to have a sizeable impact on the iteration counts despite otherwise being relatively stable to changes in $\rho$. To a lesser extent we also lose scalability with DtN($k^{4/3}$) as the iteration counts now slowly increase with $N$. On the other hand, H-GenEO($\frac{1}{2}$) has strong robustness throughout, both with respect to the wave number $k$, the contrast in the heterogeneity $\rho$, and scalability as $N$ increases. As such, we see a clear preference for H-GenEO($\frac{1}{2}$) as a stable and reliable method for heterogeneous problems.

\end{paracol}
\nointerlineskip
\begin{specialtable}[H]
	\widetable
	\centering
	\caption{Preconditioned GMRES iteration counts (above) and size of coarse space (below) for the heterogeneous increasing layers problem when using ORAS with DtN($k^{4/3}$) or H-GenEO($\frac{1}{2}$) and a varying number of subdomains. A uniform decomposition into $\sqrt{N}\times\sqrt{N}$ square subdomains is used.}
	\label{Table:DtNandH-GenEO_IncreasingLayers}
	\tabulinesep=1.2mm
	\begin{tabu}{ccc|cccccc|cccccc}
		\toprule
		& & & \multicolumn{12}{c}{Number of subdomains $N$} \\
		& & & \multicolumn{6}{c}{DtN($k^{4/3}$)} & \multicolumn{6}{c}{H-GenEO($\frac{1}{2}$)} \\
		$\omega$ & $h^{-1}$ & $\rho$ & 16 & 36 & 64 & 100 & 144 & 196 & 16 & 36 & 64 & 100 & 144 & 196 \\
		\midrule
		\multirow{2}{*}{29.3} & \multirow{2}{*}{200} & 10 & 29 & 37 & 41 & 52 & 55 & 58 & 15 & 16 & 19 & 18 & 18 & 19 \\
		& & 1000 & 44 & 44 & 50 & 58 & 52 & 52 & 15 & 15 & 17 & 18 & 17 & 17 \\
		\midrule
		\multirow{2}{*}{46.5} & \multirow{2}{*}{400} & 10 & 32 & 38 & 41 & 66 & 65 & 73 & 15 & 16 & 16 & 19 & 18 & 18 \\
		& & 1000 & 63 & 69 & 74 & 84 & 73 & 71 & 14 & 15 & 16 & 18 & 17 & 17 \\
		\midrule
		\multirow{2}{*}{73.8} & \multirow{2}{*}{800} & 10 & 35 & 43 & 42 & 40 & 58 & 69 & 15 & 17 & 16 & 17 & 18 & 17 \\
		& & 1000 & 89 & 93 & 107 & 111 & 114 & 109 & 14 & 15 & 15 & 16 & 17 & 16 \\
		\midrule\midrule
		\multirow{2}{*}{29.3} & \multirow{2}{*}{200} & 10 & 116 & 173 & 234 & 363 & 399 & 467 & 224 & 354 & 452 & 662 & 679 & 754 \\
		& & 1000 & 84 & 111 & 136 & 285 & 329 & 371 & 222 & 350 & 446 & 642 & 679 & 741 \\
		\midrule
		\multirow{2}{*}{46.5} & \multirow{2}{*}{400} & 10 & 208 & 317 & 405 & 600 & 704 & 812 & 458 & 706 & 990 & 1234 & 1523 & 1678 \\
		& & 1000 & 144 & 176 & 202 & 421 & 496 & 554 & 450 & 693 & 990 & 1216 & 1512 & 1666 \\
		\midrule
		\multirow{2}{*}{73.8} & \multirow{2}{*}{800} & 10 & 379 & 557 & 693 & 1142 & 1217 & 1404 & 930 & 1425 & 2074 & 2584 & 3060 & 3553 \\
		& & 1000 & 254 & 294 & 326 & 748 & 784 & 838 & 914 & 1409 & 2058 & 2572 & 3059 & 3534 \\
		\bottomrule
	\end{tabu}
\end{specialtable}
\begin{paracol}{2}
\switchcolumn

We also consider the case of a diagonal layers problem (Figure~\ref{Fig:DiagonalLayers}) in Table~\ref{Table:DtNandH-GenEO_DiagonalLayers}. Here the issues with DtN($k^{4/3}$) are reduced but there is still some increase in iteration counts, especially for higher wave numbers $k$ and a larger number of subdomains $N$. We note, in results not shown here, that DtN($k^{3/2}$) also suffers from the same lack of robustness. For H-GenEO($\frac{1}{2}$), however, we still have good robustness to the parameters of the problem.

We further consider the heterogeneous problem when making use of METIS for more general non-uniform subdomain decompositions in order to ensure that H-GenEO is able to handle both difficulties together. For this we consider the alternating layers problem (Figure~\ref{Fig:AlternatingLayers}) and provide results in Table~\ref{Table:H-GenEO_AlternatingLayers_METIS}. We find that H-GenEO($\frac{1}{2}$) performs very well and continues to provide a rather robust method, even in the presence of heterogeneity on non-uniform subdomains. This further evidences the strength of the H-GenEO($\frac{1}{2}$) approach and we will shortly study it more closely, dropping reference to the eigenvalue tolerance and simply denoting the method as H-GenEO. First, however, we provide a brief examination into the use of higher order discretisation.

\newpage
\end{paracol}
\nointerlineskip
\begin{specialtable}[H]
	\widetable
	\centering
	\caption{Preconditioned GMRES iteration counts (above) and size of coarse space (below) for the heterogeneous diagonal layers problem when using ORAS with DtN($k^{4/3}$) or H-GenEO($\frac{1}{2}$) and a varying number of subdomains. A uniform decomposition into $\sqrt{N}\times\sqrt{N}$ square subdomains is used.}
	\label{Table:DtNandH-GenEO_DiagonalLayers}
	\tabulinesep=1.2mm
	\begin{tabu}{ccc|cccccc|cccccc}
		\toprule
		& & & \multicolumn{12}{c}{Number of subdomains $N$} \\
		& & & \multicolumn{6}{c}{DtN($k^{4/3}$)} & \multicolumn{6}{c}{H-GenEO($\frac{1}{2}$)} \\
		$\omega$ & $h^{-1}$ & $\rho$ & 16 & 36 & 64 & 100 & 144 & 196 & 16 & 36 & 64 & 100 & 144 & 196 \\
		\midrule
		\multirow{2}{*}{29.3} & \multirow{2}{*}{200} & 10 & 13 & 14 & 13 & 14 & 21 & 25 & 16 & 18 & 20 & 18 & 23 & 25 \\
		& & 1000 & 13 & 14 & 14 & 14 & 22 & 25 & 16 & 18 & 20 & 18 & 23 & 25 \\
		\midrule
		\multirow{2}{*}{46.5} & \multirow{2}{*}{400} & 10 & 15 & 14 & 14 & 16 & 25 & 31 & 16 & 17 & 17 & 26 & 21 & 22 \\
		& & 1000 & 15 & 14 & 15 & 16 & 25 & 34 & 16 & 17 & 18 & 27 & 22 & 22 \\
		\midrule
		\multirow{2}{*}{73.8} & \multirow{2}{*}{800} & 10 & 14 & 18 & 16 & 15 & 20 & 26 & 16 & 17 & 17 & 17 & 19 & 20 \\
		& & 1000 & 15 & 18 & 16 & 15 & 32 & 39 & 16 & 17 & 17 & 17 & 19 & 20 \\
		\midrule\midrule
		\multirow{2}{*}{29.3} & \multirow{2}{*}{200} & 10 & 336 & 593 & 866 & 1090 & 1376 & 1390 & 260 & 376 & 499 & 689 & 737 & 828 \\
		& & 1000 & 336 & 594 & 866 & 1090 & 1375 & 1390 & 259 & 375 & 499 & 687 & 737 & 826 \\
		\midrule
		\multirow{2}{*}{46.5} & \multirow{2}{*}{400} & 10 & 621 & 1075 & 1540 & 1910 & 2370 & 2622 & 543 & 789 & 1095 & 1384 & 1599 & 1825 \\
		& & 1000 & 621 & 1075 & 1539 & 1907 & 2368 & 2614 & 541 & 790 & 1093 & 1381 & 1596 & 1824 \\
		\midrule
		\multirow{2}{*}{73.8} & \multirow{2}{*}{800} & 10 & 1164 & 1947 & 2692 & 3592 & 4145 & 4608 & 1145 & 1636 & 2243 & 2823 & 3233 & 3681 \\
		& & 1000 & 1163 & 1946 & 2693 & 3592 & 4131 & 4569 & 1141 & 1633 & 2239 & 2822 & 3232 & 3671 \\
		\bottomrule
	\end{tabu}
\end{specialtable}

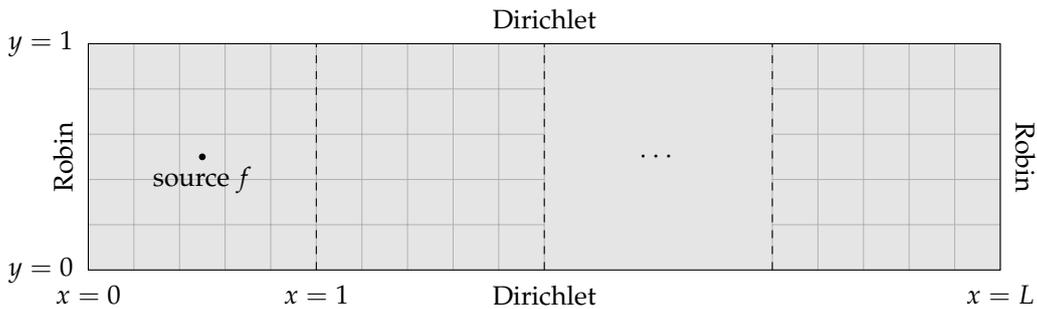
\begin{figure}[H]
	\widefigure
	\centering
	\begin{tikzpicture}[scale=1.5]
		% Grid
		\draw[step=0.4cm,gray!50,very thin, shift={(-1,-1)}] (0,0) grid (4,2);
		\draw[step=0.4cm,gray!50,very thin, shift={(5,-1)}] (0,0) grid (2,2);
		% Boundaries
		\fill[gray,opacity=0.2] (-1,1) -- (-1,-1) -- (7,-1) -- (7,1) -- cycle;
		\draw (-1,1) -- (-1,-1) -- (7,-1) -- (7,1) -- cycle;
		\draw[dashed] (1,-1) -- (1,1);
		\draw[dashed] (3,-1) -- (3,1);
		\draw[dashed] (5,-1) -- (5,1);
		% Boundary labels
		\draw (7.05,0) node[above,rotate=-90] {Robin};
		\draw (-1.05,0) node[above,rotate=90] {Robin};
		\draw (3,-1.05) node[below] {Dirichlet};
		\draw (3,1.05) node[above] {Dirichlet};
		% Coordinate labels
		\draw (-1.05,1) node[left] {$y = 1$}
		(-1.05,-1) node[left] {$y = 0$}
		(-1,-1.05) node[below] {$x = 0$}
		(1,-1.05) node[below] {$x = 1$}
		(7,-1.05) node[below] {$x = L$};
		\draw (4,0) node {$\cdots$};
		% Point source
		\draw[fill=black] (0,0) circle (0.025) node[below] {source $f$};
	\end{tikzpicture}
	\caption{Schematic of the growing 2D wave guide model problem used in a weak scaling test on $N = 25L$ fixed size subdomains, with the underlying non-overlapping subdomains shown in grey.}
	\label{Fig:WaveGuide2DWeakScaling}
\end{figure}

\begin{specialtable}[H]
	\widetable
	\centering
	\caption{Preconditioned GMRES iteration counts for the heterogeneous alternating layers problem with $\rho = 10 / 100 / 1000$ when using ORAS with H-GenEO($\frac{1}{2}$) and a varying number of subdomains. A non-uniform decomposition into $N$ subdomains is used, given by METIS.}
	\label{Table:H-GenEO_AlternatingLayers_METIS}
	\tabulinesep=1.2mm
	\begin{tabu}{cc|ccc|ccc|ccc|ccc|ccc|ccc}
		\toprule
		& & \multicolumn{18}{c}{Number of subdomains $N$ with sub-columns for $\rho = 10 / 100 / 1000$} \\
		$\omega$ & $h^{-1}$ & \multicolumn{3}{c|}{20} & \multicolumn{3}{c|}{40} & \multicolumn{3}{c|}{80} & \multicolumn{3}{c|}{120} & \multicolumn{3}{c|}{160} & \multicolumn{3}{c}{200} \\
		\midrule
		18.5 & 100 & 17 & 17 & 17 & 19 & 19 & 19 & 21 & 21 & 21 & 27 & 27 & 27 & 31 & 31 & 31 & 33 & 33 & 33 \\
		29.3 & 200 & 16 & 16 & 16 & 17 & 17 & 17 & 19 & 19 & 19 & 20 & 20 & 20 & 21 & 21 & 21 & 23 & 23 & 23 \\
		46.5 & 400 & 17 & 18 & 18 & 18 & 18 & 18 & 22 & 23 & 23 & 25 & 26 & 26 & 27 & 28 & 28 & 28 & 29 & 29 \\
		73.8 & 800 & 16 & 16 & 16 & 17 & 17 & 17 & 18 & 18 & 18 & 18 & 19 & 19 & 19 & 20 & 20 & 23 & 23 & 23 \\
		117.2 & 1600 & 15 & 15 & 15 & 15 & 16 & 16 & 16 & 16 & 16 & 16 & 16 & 16 & 16 & 16 & 16 & 16 & 16 & 16 \\
		\bottomrule
	\end{tabu}
\end{specialtable}
\begin{paracol}{2}
\switchcolumn

\subsection{Higher Order Finite Elements}

We briefly investigate the use of higher order finite elements. In particular, we consider the use of P2 elements as opposed to P1 elements. To give a direct comparison we utilise the same meshes and in Table~\ref{Table:DtNandH-GenEO_DiagonalLayers_P2} give results for the heterogeneous diagonal layers problem with $\rho = 10$; equivalent results for P1 elements are given in Table~\ref{Table:DtNandH-GenEO_DiagonalLayers}. We observe that both iteration counts and coarse space sizes remain rather similar to the case of P1 elements. This suggests that, by itself, the order of the underlying finite elements used does not strongly affect performance. However, the typical meshes employed for higher order elements may be coarser and further studies would be required to observe how this affects the utility of the spectral coarse spaces presented here. Further investigation into which method and order of finite element approximation provides the most efficient choice is beyond the scope of the present study.

\newpage
\end{paracol}
\nointerlineskip
\begin{specialtable}[H]
\widetable
\centering
\caption{Preconditioned GMRES iteration counts (above) and size of coarse space (below) for P2 finite element discretisation of the heterogeneous diagonal layers problem with $\rho = 10$ when using ORAS with DtN($k^{4/3}$) or H-GenEO($\frac{1}{2}$) and a varying number of subdomains. A uniform decomposition into $\sqrt{N}\times\sqrt{N}$ square subdomains is used.}
\label{Table:DtNandH-GenEO_DiagonalLayers_P2}
\tabulinesep=1.2mm
\begin{tabu}{ccc|cccccc|cccccc}
	\toprule
	& & & \multicolumn{12}{c}{Number of subdomains $N$} \\
	& & & \multicolumn{6}{c}{DtN($k^{4/3}$)} & \multicolumn{6}{c}{H-GenEO($\frac{1}{2}$)} \\
	$\omega$ & $h^{-1}$ & $\rho$ & 16 & 36 & 64 & 100 & 144 & 196 & 16 & 36 & 64 & 100 & 144 & 196 \\
	\midrule
	18.5 & 100 & 10 & 13 & 11 & 10 & 10 & 10 & 18 & 15 & 16 & 19 & 18 & 23 & 24 \\
	29.3 & 200 & 10 & 14 & 12 & 12 & 13 & 19 & 25 & 15 & 17 & 18 & 18 & 23 & 25 \\
	46.5 & 400 & 10 & 15 & 12 & 12 & 15 & 23 & 30 & 15 & 16 & 17 & 20 & 21 & 22 \\
	73.8 & 800 & 10 & 17 & 16 & 14 & 13 & 18 & 25 & 15 & 16 & 16 & 17 & 18 & 20 \\
	\midrule\midrule
	18.5 & 100 & 10 & 151 & 326 & 510 & 706 & 898 & 937 & 125 & 186 & 231 & 346 & 398 & 519 \\
	29.3 & 200 & 10 & 300 & 608 & 892 & 1086 & 1444 & 1516 & 260 & 377 & 507 & 686 & 733 & 824 \\
	46.5 & 400 & 10 & 589 & 1108 & 1572 & 2069 & 2491 & 2638 & 540 & 794 & 1100 & 1403 & 1594 & 1820 \\
	73.8 & 800 & 10 & 919 & 1916 & 2862 & 3614 & 4409 & 4748 & 1144 & 1645 & 2239 & 2805 & 3225 & 3695 \\
	\bottomrule
\end{tabu}
\end{specialtable}
\begin{paracol}{2}
\switchcolumn

\subsection{The Effect of Boundary Conditions within the H-GenEO Eigenproblem}

We now consider the choice of boundary conditions within H-GenEO in light of the fact that, for wave propagation problems, impedance conditions can often prove more practical within overlapping Schwarz methods. To this end, we consider the H-GenEO eigenproblem where the Neumann boundary condition is replaced by an impedance condition instead (i.e., the Robin condition in \eqref{ORASLocalRobinBC}). Results for this impedance-H-GenEO method are given in Table~\ref{Table:ImpedanceH-GenEO} for the homogeneous problem with uniform square subdomains. We see that, while the use of this eigenproblem retains the good behaviour of H-GenEO as $k$ increases, it lacks scalability as we increase the number of subdomains $N$. We note that the size of the coarse space is very similar to that when the standard Neumann condition is used (a direct comparison can be made for $k=73.8$ with results in Table~\ref{Table:DtNandH-GenEO_Scalability}) and so this is not simply an artefact of a smaller coarse space. This shows that the Neumann condition within the eigenproblem is an important aspect of H-GenEO.

\end{paracol}
\nointerlineskip
\begin{specialtable}[H]
	\widetable
	\centering
	\caption{Preconditioned GMRES iteration counts (above) and size of coarse space (below) for the homogeneous problem when using ORAS with impedance-H-GenEO (the eigenproblem \eqref{DiscreteH-GenEOEigenproblem} is altered to have impedance as opposed to Neumann boundary conditions on the left-hand side) and a varying number of subdomains. A uniform decomposition into $\sqrt{N}\times\sqrt{N}$ square subdomains is used.}
	\label{Table:ImpedanceH-GenEO}
	\tabulinesep=1.2mm
	\begin{tabu}{cc|ccccccccccccc}
		\toprule
		& & \multicolumn{13}{c}{Number of subdomains $N$} \\
		$k$ & $h^{-1}$ & 4 & 9 & 16 & 25 & 36 & 49 & 64 & 81 & 100 & 121 & 144 & 169 & 196 \\
		\midrule
		18.5 & 100 & 17 & 19 & 23 & 27 & 30 & 36 & 42 & 45 & 43 & 58 & 61 & 61 & 67 \\
		29.3 & 200 & 17 & 19 & 22 & 25 & 34 & 33 & 41 & 38 & 35 & 49 & 60 & 62 & 65 \\
		46.5 & 400 & 15 & 18 & 19 & 22 & 26 & 25 & 26 & 39 & 43 & 47 & 52 & 51 & 54 \\
		73.8 & 800 & 15 & 19 & 19 & 20 & 25 & 27 & 26 & 35 & 33 & 39 & 43 & 44 & 51 \\
		\midrule\midrule
		18.5 & 100 & 68 & 102 & 140 & 158 & 204 & 217 & 236 & 287 & 344 & 341 & 404 & 477 & 504 \\
		29.3 & 200 & 148 & 215 & 296 & 370 & 392 & 521 & 504 & 576 & 720 & 768 & 725 & 793 & 908 \\
		46.5 & 400 & 360 & 492 & 628 & 754 & 917 & 988 & 1236 & 1176 & 1468 & 1550 & 1740 & 1807 & 1930 \\
		73.8 & 800 & 848 & 1106 & 1420 & 1696 & 1877 & 2218 & 2432 & 2574 & 2960 & 3180 & 3443 & 3834 & 3732 \\
		\bottomrule
	\end{tabu}
\end{specialtable}
\begin{paracol}{2}
\switchcolumn

\subsection{The Effect of More Overlap when Using H-GenEO}

So far all our results use minimal overlap. Here we consider the case of increasing the overlap between subdomains, this is done by adding on layers of adjoining elements to each subdomain in a symmetric way, so that minimal overlap is given by an overlap parameter of 2; that is the overlapping region has a width of 2 elements. In Table~\ref{Table:H-GenEO_Overlap} we report results for increasing overlap when using H-GenEO for the homogeneous problem with $k = 46.5$ and $h^{-1} = 400$ and where a uniform decomposition is used. We see that adding on a small amount of overlap can slightly decrease the iteration counts but increasing the overlap width further can give much poorer results, especially when using a large number of subdomains $N$. One possible explanation could be an increase in the ``colouring constant'' (see, e.g., \cite[Definition 5.5]{Dolean:15:DDM}). We note that the size of the coarse space decreases somewhat as the overlap is increased, however, the extra computational effort required to deal with the larger subdomains will hamper any gains from this, along with the increased iteration counts. From these results we determine that the H-GenEO coarse space is best suited to the case of minimal overlap, as we have used elsewhere throughout this work.

\end{paracol}
\nointerlineskip
\begin{specialtable}[H]
	\widetable
	\centering
	\caption{Preconditioned GMRES iteration counts (above) and size of coarse space (below) for the homogeneous problem when using ORAS with H-GenEO, varying the amount of overlap (in terms of element width, with 2 representing minimal overlap) and number of subdomains for $k = 46.5$ and $h^{-1} = 400$. A uniform decomposition into $\sqrt{N}\times\sqrt{N}$ square subdomains is used.}
	\label{Table:H-GenEO_Overlap}
	\tabulinesep=1.2mm
	\begin{tabu}{c|ccccccccccccc}
		\toprule
		& \multicolumn{13}{c}{Number of subdomains $N$} \\
		Overlap & 4 & 9 & 16 & 25 & 36 & 49 & 64 & 81 & 100 & 121 & 144 & 169 & 196 \\
		\midrule
		2 & 14 & 15 & 15 & 17 & 16 & 16 & 16 & 20 & 26 & 19 & 19 & 22 & 21 \\
		4 & 10 & 11 & 11 & 12 & 12 & 13 & 12 & 17 & 14 & 16 & 17 & 21 & 20 \\
		8 & 8 & 10 & 10 & 10 & 13 & 13 & 12 & 20 & 23 & 22 & 26 & 31 & 27 \\
		16 & 13 & 21 & 26 & 26 & 37 & 77 & 61 & 75 & 86 & 109 & 178 & 157 & 164 \\
		\midrule\midrule
		2 & 368 & 492 & 644 & 779 & 938 & 1030 & 1248 & 1195 & 1476 & 1558 & 1758 & 1845 & 2016 \\
		4 & 352 & 472 & 600 & 699 & 871 & 947 & 1088 & 1124 & 1296 & 1449 & 1689 & 1697 & 1690 \\
		8 & 336 & 436 & 538 & 650 & 799 & 863 & 1024 & 981 & 1132 & 1395 & 1511 & 1425 & 1512 \\
		16 & 316 & 417 & 500 & 610 & 733 & 717 & 920 & 942 & 1108 & 1239 & 1086 & 1212 & 1280 \\
		\bottomrule
	\end{tabu}
\end{specialtable}
\begin{paracol}{2}
\switchcolumn

\subsection{Weak Scalability and Timing Results for H-GenEO}

In this section we consider weak scalability of the H-GenEO method. To approach this we consider a growing wave guide domain where fixed size subdomains are added, each with the same number of degrees of freedom (dofs). This is done by repeatedly adding a unit square, which is split into 25 non-overlapping square subdomains, to the right of the existing domain $L$ times to give $\Omega = (0,L)\times(0,1)$. Along the long edges of the global domain we prescribe homogeneous Dirichlet boundary conditions while Robin conditions are used at each end of the wave guide; a schematic for this weak scaling test is given in Figure~\ref{Fig:WaveGuide2DWeakScaling}. Heterogeneity is given by the alternating layers problem (see Figure~\ref{Fig:AlternatingLayers}) with the ten layers extending across the length of the wave guide.

To deal with the large problem sizes (reaching up to $10,\!253,\!601$ dofs) and provide appropriate timing results, we assign one core per subdomain and solve using the ARCHIE-WeSt supercomputing facility on up to 400 cores (the machine uses Intel Xeon Gold 6138 processors at $2.0\,$GHz with $4.8\,$GB RAM per core).

\newpage
\end{paracol}
\nointerlineskip
\begin{paracol}{2}
\switchcolumn

In Table~\ref{Table:H-GenEO_Weak_Scalability} we give results for a length of domain $L = 2$ up to $L = 16$, that is from $N = 50$ to $N = 400$ subdomains, using $k = 73.8$, $h^{-1} = 800$ and $\rho = 100$. We observe that the iteration counts remain almost constant, increasing only very mildly with $N$. On the other hand, the coarse space size grows linearly with $N$, as expected given subdomains have a fixed size, and we see that the time taken to solve the eigenproblems stays constant. However, due the increasing size of the coarse space, the run times slowly increase so that the efficiency degrades as we solve larger and larger problems. Further, the setup of the decomposition and partitioning is performed sequentially and so this setup time increases with $N$. Nonetheless, when removing this initial setup the efficiency still slowly reduces as we increase the problem size, as can be seen in the final row of Table~\ref{Table:H-GenEO_Weak_Scalability}. To overcome this, the factorisation of the coarse space operator and associated solves must be made scalable. While a multi-level approach is therefore an attractive option, it is not yet clear how to formulate such a strategy for Helmholtz problems. Additionally, the coarse problem is only assigned to one process here and so spreading it over more processes as $N$ increases may ease the coarse solve becoming a bottleneck.

\end{paracol}
\nointerlineskip
\begin{specialtable}[H]
	\widetable
	\centering
	\caption{Weak scaling results and timings for the alternating layers problem when using ORAS with H-GenEO($\frac{1}{2}$) and a varying number of subdomains for $k = 73.8$, $h^{-1} = 800$ and $\rho = 100$. A uniform decomposition into $N = 25L$ subdomains is used, as depicted in Figure~\ref{Fig:WaveGuide2DWeakScaling}. Note that setup refers to the initial decomposition and partitioning, which is performed sequentially, while the local problems and eigensolves are carried out in parallel.}
	\label{Table:H-GenEO_Weak_Scalability}
	\tabulinesep=1.2mm
	\begin{tabu}{c|cccccccc}
		\toprule
		$N$ & 50 & 100 & 150 & 200 & 250 & 300 & 350 & 400 \\
		\midrule
		Iteration count & 17 & 18 & 18 & 19 & 19 & 20 & 21 & 21 \\
		Coarse space size & 3010 & 6150 & 9290 & 12430 & 15570 & 18710 & 21850 & 24990 \\
		Total run time (s) & 45.8 & 48.6 & 53.0 & 58.7 & 63.5 & 70.0 & 79.7 & 88.1 \\
		Weak scaling efficiency & $-$ & 94.2\% & 86.4\% & 78.0\% & 72.1\% & 65.4\% & 57.5\% & 52.0\% \\
		Eigensolve time (s) & 37.1 & 37.9 & 37.9 & 38.3 & 37.8 & 37.9 & 37.9 & 37.7 \\
		Setup time (s) & 5.5 & 7.7 & 12.9 & 16.5 & 19.9 & 23.8 & 27.4 & 30.8 \\
		Efficiency without setup & $-$ & 98.5\% & 100.5\% & 95.5\% & 92.4\% & 87.2\% & 77.1\% & 70.3\% \\
		\bottomrule
	\end{tabu}
\end{specialtable}
\begin{paracol}{2}
\switchcolumn

\subsection{High Wave Number Strong Scalability and Timing Results for H-GenEO}

To conclude our numerical results we consider the use of H-GenEO within a high wave number example and explore timings and the overall strong scalability of the approach. To this end we consider the original homogeneous wave guide, as outlined in Figure~\ref{Fig:WaveGuide2D}, with a wave number of $k = 186.0$ and mesh width $h^{-1} = 3200$, giving a total problem size of $10,\!246,\!401$ dofs. To solve we use METIS to give a non-uniform decomposition into subdomains and again utilise the ARCHIE-WeSt supercomputing facility to assign one core per subdomain.

Results are tabulated in Table~\ref{Table:H-GenEO_Scalability_Timings_METIS}, in which we detail the run time (in seconds) and the percentage of time spent by the eigensolver and coarse factorisation. The run times are further used to determine the parallel efficiency based on the smallest run on $N = 80$ cores. Timing results are also displayed graphically in Figure~\ref{Fig:Timings}, where we that the bulk of the computation is spent in the eigensolves and factorisation of the coarse problem; as such, while the METIS decomposition is sequential, we do not separate out the setup time here given that it is comparatively rather small. From Table~\ref{Table:H-GenEO_Scalability_Timings_METIS} we see that the iteration counts show good scalability, increasing only very mildly as we increase the number of subdomains fivefold. This is also seen in the run times, which decrease as we use more cores, and hence more subdomains, to solve the problem. In particular, we see that the parallel efficiency remains over 100\%, showing strong scalability of the approach, in part due to the fact that as $N$ increases we have to solve smaller eigenproblems and so the percentage of time spent by the eigensolver drops significantly as we increase $N$. Overall these provide promising results that the H-GenEO method can be effective for the solution of high wave number problems in 2D.

\end{paracol}
\nointerlineskip
\begin{specialtable}[H]
	\widetable
	\centering
	\caption{Strong scaling results and timings for the homogeneous problem when using ORAS with H-GenEO($\frac{1}{2}$) and a varying number of subdomains for $k = 186.0$ and $h^{-1} = 3200$, giving a total of $10,\!246,\!401$ dofs. A non-uniform decomposition into $N$ subdomains is used, given by METIS. The average local eigenproblem size is given approximately as the number of dofs divided by $N$.}
	\label{Table:H-GenEO_Scalability_Timings_METIS}
	\tabulinesep=1.2mm
	\begin{tabu}{c|ccccccccc}
		\toprule
		$N$ & 80 & 120 & 160 & 200 & 240 & 280 & 320 & 360 & 400 \\
		\midrule
		Iteration count & 14 & 16 & 15 & 16 & 17 & 17 & 18 & 19 & 19 \\
		Coarse space size & 16014 & 19018 & 21348 & 23747 & 25560 & 27270 & 28793 & 30357 & 31773 \\
		Total run time (s) & 1214.4 & 614.6 & 404.4 & 279.3 & 217.3 & 195.0 & 159.4 & 154.0 & 147.6 \\
		Parallel efficiency & $-$ & 132\% & 150\% & 174\% & 186\% & 178\% & 190\% & 175\% & 165\% \\
		\midrule
		Eigenproblem size (approx.) & 128080 & 85387 & 64040 & 51232 & 42693 & 36594 & 32020 & 28462 & 25616 \\
		Average no. of eigenvectors & 200.2 & 158.5 & 133.4 & 118.7 & 106.5 & 97.4 & 90.0 & 84.3 & 79.4 \\
		Eigensolve time & 68.6\% & 70.1\% & 66.8\% & 62.7\% & 53.7\% & 52.3\% & 41.6\% & 39.3\% & 37.3\% \\
		Coarse factorisation time & 29.6\% & 27.8\% & 30.8\% & 34.4\% & 42.8\% & 44.3\% & 54.0\% & 56.3\% & 58.3\% \\
		\bottomrule
	\end{tabu}
\end{specialtable}

\begin{figure}[H]
	\widefigure
	\centering
	\includegraphics[width=0.48\textwidth,trim=0cm 0cm 0cm 0cm ,clip]{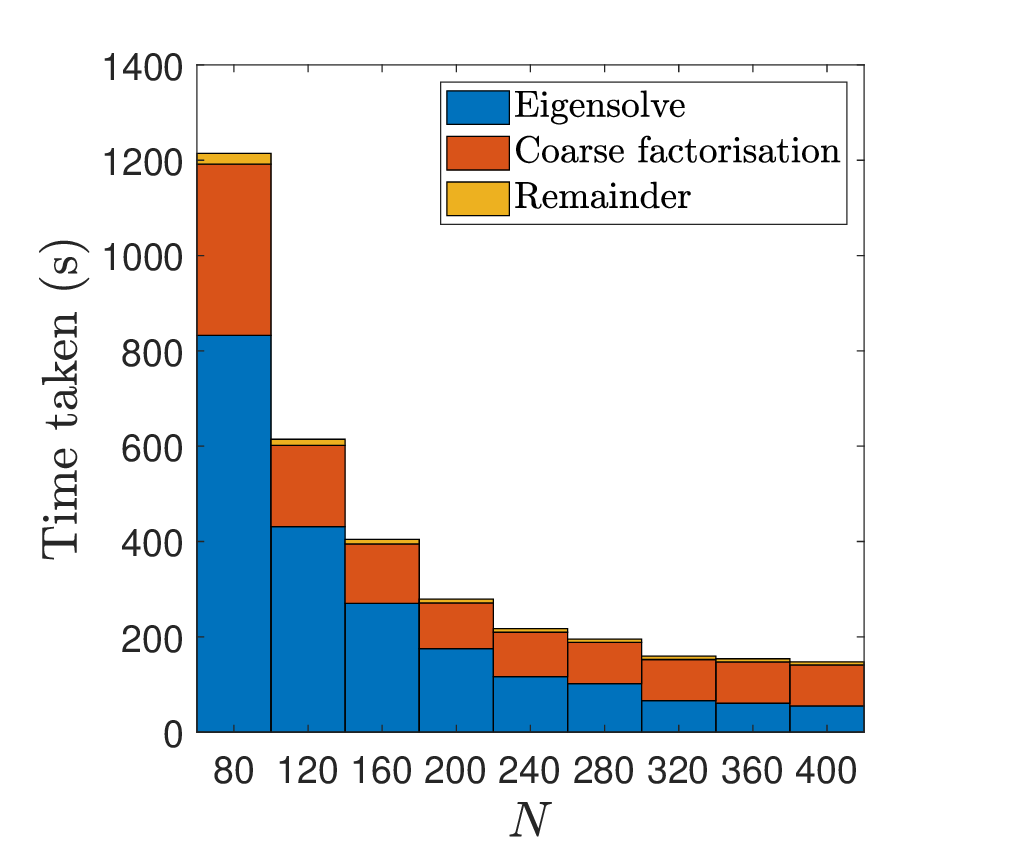} \includegraphics[width=0.48\textwidth,trim=0cm 0cm 0cm 0cm ,clip]{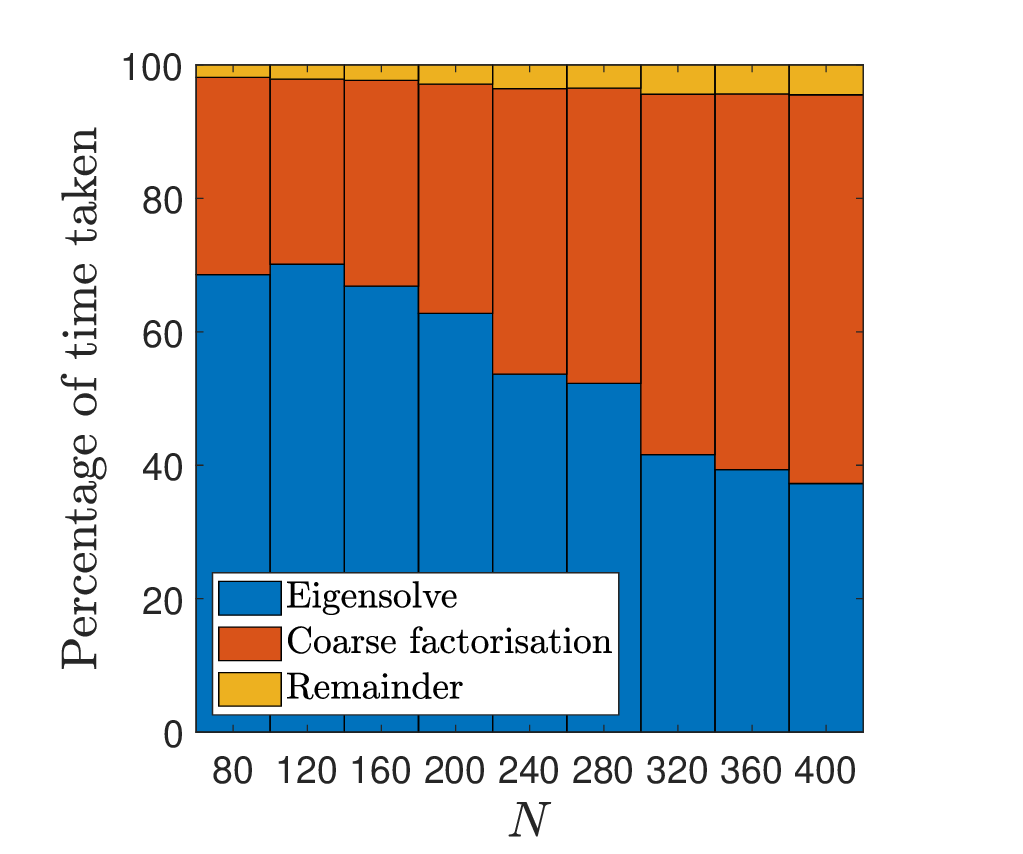}
	\caption{Timings for the homogeneous problem when using ORAS with H-GenEO($\frac{1}{2}$) and a varying number of subdomains for $k = 186.0$ and $h^{-1} = 3200$, giving a total of $10,\!246,\!401$ dofs. A non-uniform decomposition into $N$ subdomains is used, given by METIS.}
	\label{Fig:Timings}
\end{figure}
\begin{paracol}{2}
\switchcolumn

%%%%%%%%%%%%%%%%%%%%%%%%%%%%%%%%%%%%%%%%%%
\section{Conclusions}
\label{sec:Conclusions}

In this work we have developed and explored a GenEO-type coarse space for additive Schwarz methods that is appropriate for the heterogeneous Helmholtz problem. We have conducted extensive numerical tests to show how this approach behaves on a 2D model test problem of a wave guide discretised using finite elements on a pollution-free mesh, comparing our method with the DtN coarse space. We find that only our H-GenEO approach is robust to heterogeneity and increasing wave number $k$, and further provides good scaling such that iteration counts of right-preconditioned GMRES are only mildly dependent on the number of subdomains. This dependence is strongest for lower wave numbers on many subdomains and reduces as $k$ grows. Furthermore, convergence does not deteriorate with increasing wave number, albeit at the cost of a coarse space which grows as $k$ increases. This is achieved consistently for non-uniform partitioning into subdomains and in the presence of strong heterogeneity. These results show promise that H-GenEO can be used as an effective coarse space for challenging heterogeneous problems.

Finally, we discuss how these findings differ from that in the companion paper \cite{Bootland:2021:ACO}, where none of the approaches (including the spectral coarse spaces explored in detail here) as seen to be clearly favourable over a wide range of problem settings. A key difference is that well-resolved pollution-free meshes are used in the present work, while the more realistic benchmark problems examined in \cite{Bootland:2021:ACO} use under-resolved meshes with a fixed number of points per wavelength, as is typical in engineering practice. As well the more complex test cases considered, another contributing factor to the difference in studies is that in \cite{Bootland:2021:ACO} a fixed number of eigenvectors was taken per subdomain (limiting the size of the coarse space), whereas here eigenvalue thresholding was employed to provide a more robust approach, at the cost of a potentially larger coarse space. As such, the present study can be thought of as investigating the extent of what can be achieved in terms of a robust method in the ideal case, while \cite{Bootland:2021:ACO} presents a viewpoint within the regime of more challenging practical problems. Identifying coarse spaces which can bridge this gap and provide a consistently robust approach even for the most demanding real-world applications remains an open area of research.

%%%%%%%%%%%%%%%%%%%%%%%%%%%%%%%%%%%%%%%%%%
\vspace{6pt}

%%%%%%%%%%%%%%%%%%%%%%%%%%%%%%%%%%%%%%%%%%
%% optional
%\supplementary{The following are available online at \linksupplementary{s1}, Figure S1: title, Table S1: title, Video S1: title.}

% Only for the journal Methods and Protocols:
% If you wish to submit a video article, please do so with any other supplementary material.
% \supplementary{The following are available at \linksupplementary{s1}, Figure S1: title, Table S1: title, Video S1: title. A supporting video article is available at doi: link.}

%%%%%%%%%%%%%%%%%%%%%%%%%%%%%%%%%%%%%%%%%%
\authorcontributions{Conceptualization, N.B. and V.D.; methodology, N.B.; software, N.B.; validation, N.B.; formal analysis, N.B. and V.D.; investigation, N.B.; resources, V.D.; writing---original draft preparation, N.B.; writing---review and editing, N.B. and V.D.; visualization, N.B.; supervision, V.D.; project administration, V.D.; funding acquisition, V.D. All authors have read and agreed to the published version of the manuscript.}

\funding{This research was funded by EPSRC grant number EP/S004017/1.}

\acknowledgments{Numerical results were obtained using the ARCHIE-WeSt High Performance Computer (\href{www.archie-west.ac.uk}{www.archie-west.ac.uk}) based at the University of Strathclyde.}

\conflictsofinterest{The authors declare no conflict of interest.}

%%%%%%%%%%%%%%%%%%%%%%%%%%%%%%%%%%%%%%%%%%
%% Only for journal Encyclopedia
%\entrylink{The Link to this entry published on the encyclopedia platform.}

%%%%%%%%%%%%%%%%%%%%%%%%%%%%%%%%%%%%%%%%%%
%% Optional
%\abbreviations{Abbreviations}{The following abbreviations are used in this manuscript:\\
%
%\noindent
%\begin{tabular}{@{}ll}
%MDPI & Multidisciplinary Digital Publishing Institute\\
%DOAJ & Directory of open access journals\\
%TLA & Three letter acronym\\
%LD & Linear dichroism
%\end{tabular}}

\end{paracol}
\reftitle{References}


\begin{thebibliography}{999}

\bibitem{Amestoy:2001:AFA}
Amestoy, P.~R.; Duff, I.~S.; L'Excellent, J.-Y.; Koster, J. A fully asynchronous multifrontal solver using distributed dynamic scheduling. {\em SIAM J. Matrix Anal. Appl.} {\bf 2001}, {\em 23}, 15--41.

\bibitem{Babuska:1997:IPE}
Babuska, I.~M.; Sauter, S.~A. Is the pollution effect of the {FEM} avoidable for the {Helmholtz} equation considering high wave numbers? {\em SIAM J. Numer. Anal.} {\bf 1997}, {\em 34}, 2392--2423.

\bibitem{Beriot:2021:AAP}
Beriot, H.; Modave, A. An automatic perfectly matched layer for acoustic finite element simulations in convex domains of general shape. {\em Internat. J. Numer. Methods Engrg.} {\bf 2021}, {\em 122}, 1239--1261.

\bibitem{Bonazzoli:2019:AOP}
Bonazzoli, M.; Dolean, V.; Graham, I.~G.; Spence, E.~A.; Tournier, P.-H. Domain decomposition preconditioning for the high-frequency time-harmonic {Maxwell} equations with absorption. {\em Math. Comp.} {\bf 2019}, {\em 88}, 2559--2604.

\bibitem{Bootland:2019:ODN}
Bootland, N.; Dolean, V. On the {Dirichlet-to-Neumann} coarse space for solving the {Helmholtz} problem using domain decomposition. In {\em Numerical Mathematics and Advanced Applications ENUMATH 2019}; Vermolen, F.~J., Vuik, C., Eds.; Springer: Cham, Switzerland, 2021; pp. 175--184.

\bibitem{Bootland:2022:GCS}
Bootland, N.; Dolean, V.; Graham, I.~G.; Ma, C.; Scheichl, R. {GenEO} coarse spaces for heterogeneous indefinite elliptic problems. In {\em Domain Decomposition Methods in Science and Engineering XXVI}; Brenner, S., et al., Eds.; Springer: Cham, Switzerland (accepted).

\bibitem{Bootland:2021:OSM}
Bootland, N.; Dolean, V.; Graham, I.~G.; Ma, C.; Scheichl, R. Overlapping {Schwarz} methods with {GenEO} coarse spaces for indefinite and non-self-adjoint problems. \textit{arXiv preprint} {\bf 2021}, arXiv:2110.13537.

\bibitem{Bootland:2021:ACO}
Bootland, N.; Dolean, V.; Jolivet, P.; Tournier, P.-H. A comparison of coarse spaces for {Helmholtz} problems in the high frequency regime. {\em Comput. Math. Appl.} {\bf 2021}, {\em 98}, 239--253.

\bibitem{Bootland:2022:AOP}
Bootland, N.; Dolean, V.; Kyriakis, A.; Pestana, J. Analysis of parallel {Schwarz} algorithms for time-harmonic problems using block {Toeplitz} matrices. {\em Electron. Trans. Numer. Anal.} {\bf 2022}, {\em 55}, 112--141.

\bibitem{Boubendir:2012:AQO}
Boubendir, Y.; Antoine, X.; Geuzaine, C. A quasi-optimal non-overlapping domain decomposition algorithm for the {Helmholtz} equation. {\em J. Comput. Phys.} {\bf 2012}, {\em 231}, 262--280.

\bibitem{Cai:1998:OSA}
Cai, X.-C.; Casarin, M.~A.; Elliott Jr, F.~W.; Widlund, O.~B. Overlapping Schwarz algorithms for solving Helmholtz's equation. {\em Contemp. Math.} {\bf 1998}, {\em 218}, 391--399.

\bibitem{Calandra:2013:AIT}
Calandra, H.; Gratton, S.; Pinel, X.; Vasseur, X. An improved two-grid preconditioner for the solution of three-dimensional {Helmholtz} problems in heterogeneous media. {\em Numer. Linear Algebra Appl.} {\bf 2013}, {\em 20}, 663--688.

\bibitem{Claeys:2020:RTO}
Claeys, X.; Parolin, E. Robust treatment of cross points in {Optimized Schwarz Methods}. \textit{arXiv preprint} {\bf 2020}, arXiv:2003.06657.

\bibitem{Cocquet:2017:HLA}
Cocquet, P.-H.; Gander, M.~J. How large a shift is needed in the shifted {Helmholtz} preconditioner for its effective inversion by multigrid? {\em SIAM J. Sci. Comput.} {\bf 2017}, {\em 39}, A438--A478.

\bibitem{Collino:2000:DDM}
Collino, F.; Ghanemi, S.; Joly, P. Domain decomposition method for harmonic wave propagation: a general presentation. {\em Comput. Methods Appl. Mech. Engrg.} {\bf 2000}, {\em 184}, 171--211.

\bibitem{Conen:2014:ACS}
Conen, L.; Dolean, V.; Krause, R.; Nataf, F. A coarse space for heterogeneous {Helmholtz} problems based on the {Dirichlet-to-Neumann} operator. {\em J. Comput. Appl. Math.} {\bf 2014}, {\em 271}, 83--99.

\bibitem{Dai:2022:MSP}
Dai, R.; Modave, A.; Remacle, J.-F.; Geuzaine, C. Multidirectional sweeping preconditioners with non-overlapping checkerboard domain decomposition for {Helmholtz} problems. {\em J. Comput. Phys.} {\bf 2022}, {\em 453}, 110887.

\bibitem{Despres:1990:DDM}
Despr{\'e}s, B. Domain decomposition method for the {Helmholtz} problem. {\em C. R. Math. Acad. Sci. Paris. I Math.} {\bf 1990}, {\em 311}, 313--316.

\bibitem{Dolean:15:DDM}
Dolean, V.; Jolivet, P.; Nataf, F. \textit{An Introduction to Domain Decomposition Methods: Algorithms, Theory, and Parallel Implementation}; Society for Industrial and Applied Mathematics (SIAM): Philadelphia, PA, USA, 2015.

\bibitem{Dwarka:2020:HLA}
Dwarka, V.; Vuik, C. Scalable convergence using two-level deflation preconditioning for the {Helmholtz} equation. {\em SIAM J. Sci. Comput.} {\bf 2020}, {\em 42}, A901--A928.

\bibitem{Engquist:2011:SP1}
Engquist, B.; Ying, L. Sweeping preconditioner for the {Helmholtz} equation: hierarchical matrix representation. {\em Comm. Pure Appl. Math.} {\bf 2011}, {\em 64}, 697--735.

\bibitem{Engquist:2011:SP2}
Engquist, B.; Ying, L. Sweeping preconditioner for the {Helmholtz} equation: moving perfectly matched layers. {\em Multiscale. Model. Simul.} {\bf 2011}, {\em 9}, 686--710.

\bibitem{Erlangga:2004:OAC}
Erlangga, Y.; Vuik, C.; Oosterlee, C.~W. On a class of preconditioners for solving the {Helmholtz} equation. {\em Appl. Numer. Math.} {\bf 2004}, {\em 50}, 409--425.

\bibitem{Erlangga:2006:ANM}
Erlangga, Y.; Oosterlee, C.~W.; Vuik, C. A novel multigrid based preconditioner for heterogeneous {Helmholtz} problems. {\em SIAM J. Sci. Comput.} {\bf 2006}, {\em 27}, 1471--1492.

\bibitem{Ernst:2012:NAM}
Ernst, O.~G.; Gander, M.~J. Why it is difficult to solve {Helmholtz} problems with classical iterative methods. In {\em Numerical Analysis of Multiscale Problems}; Graham, I.~G., Hou, T.~Y., Lakkis, O., Scheichl, R., Eds.; Springer: Berlin, Germany, 2012; pp. 325--363.

\bibitem{Farhat:2005:FET}
Farhat, C.; Avery, P.; Tezaur, R.; Li, J. {FETI-DPH}: a dual-primal domain decomposition method for acoustic scattering. {\em J. Comput. Acoust.} {\bf 2005}, {\em 13}, 499--524.

\bibitem{Farhat:2000:ATL}
Farhat, C.; Macedo, A.; Lesoinne, M. A two-level domain decomposition method for the iterative solution of high frequency exterior {Helmholtz} problems. {\em Numer. Math.} {\bf 2000}, {\em 85}, 283--308.

\bibitem{Fish:2000:GBT}
Fish, J.; Qu, Y. Global-basis two-level method for indefinite systems. {Part 1:	convergence studies}. {\em Internat. J. Numer. Methods Engrg.} {\bf 2000}, {\em 49}, 439--460.

\bibitem{Gander:2002:OSM}
Gander, M.~J.; Magoules, F.; Nataf, F. Optimized {Schwarz} methods without overlap for the {Helmholtz} equation. {\em SIAM J. Sci. Comput.} {\bf 2002}, {\em 24}, 38--60.

\bibitem{Gander:2019:SIREV}
Gander, M.~J.; Zhang, H. A class of iterative solvers for the {Helmholtz} equation: Factorizations, sweeping preconditioners, source transfer, single layer potentials, polarized traces, and optimized {Schwarz} methods. {\em SIAM Review} {\bf 2019}, {\em 60}, 3--76.

\bibitem{Gander:2016:OSM}
Gander, M.~J.; Zhang, H. Optimized Schwarz methods with overlap for the Helmholtz equation. {\em SIAM J. Sci. Comput.} {\bf 2016}, {\em 38}, A3195--A3219.

\bibitem{Gillman:2015:BIT}
Gillman, A.; Barnett, A.~H.; Martinsson, P.-G. A spectrally accurate direct solution technique for frequency-domain scattering problems with variable media. {\em BIT} {\bf 2015}, {\em 55}, 141--170.

\bibitem{Gong:2021:COP}
Gong, S.; Gander, M.~J.; Graham, I.~G.; Lafontaine, D.; Spence, E.~A. Convergence of parallel overlapping domain decomposition methods for the {Helmholtz} equation. \textit{arXiv preprint} {\bf 2021}, arXiv:2106.05218.

\bibitem{Gong:2021:DDP}
Gong, S.; Graham, I.~G.; Spence, E.~A. Domain decomposition preconditioners for high-order discretizations of the heterogeneous {Helmholtz} equation. {\em IMA J. Numer. Anal.} {\bf 2021}, {\em 41}, 2139--2185.

\bibitem{Graham:2017:DDP}
Graham, I.~G.; Spence, E.~A.; Vainikko, E. Domain decomposition preconditioning for high-frequency {Helmholtz} problems with absorption. {\em Math. Comp.} {\bf 2017}, {\em 86}, 2089--2127.

\bibitem{Graham:2017:RRO}
Graham, I.~G.; Spence, E.~A.; Vainikko, E. Recent results on domain decomposition preconditioning for the high-frequency {Helmholtz} equation using absorption. In {\em Modern Solvers for Helmholtz Problems}; Lahaye, D., Tang, J., Vuik, K., Eds.; Birkh\"{a}user: Cham, Switzerland, 2017; pp. 3--26.

\bibitem{Graham:2020:DDW}
Graham, I.~G.; Spence, E.~A.; Zou, J. Domain Decomposition with local impedance conditions for the {Helmholtz} equation with absorption. {\em SIAM J. Numer. Anal.} {\bf 2020}, {\em 58}, 2515--2543.

\bibitem{Haferssas:2017:AAS}
Haferssas, R.; Jolivet, P.; Nataf, F. An additive Schwarz method type theory for {Lions's} algorithm and a symmetrized optimized restricted additive Schwarz method. {\em SIAM J. Sci. Comput.} {\bf 2017}, {\em 39}, A1345--A1365.

\bibitem{Harari:2000:AAN}
Harari, I.; Slavutin, M.; Turkel, E. Analytical and numerical studies of a finite element {PML} for the {Helmholtz} equation. {\em J. Comput. Acoust.} {\bf 2000}, {\em 8}, 121--137.

\bibitem{Hecht:2012:NDI}
Hecht, F. New development in {FreeFem++}, {\em J. Numer. Math.} {\bf 2012}, {\em 20}, 251--266.

\bibitem{Hu:2016:SPD}
Hu, Q.; Zhang, H. Substructuring preconditioners for the systems arising from plane wave discretization of {Helmholtz} equations. {\em SIAM J. Sci. Comput.} {\bf 2016}, {\em 38}, A2232--A2261.

\bibitem{Karypis:1998:AFA}
Karypis, G.; Kumar, V. A fast and high quality multilevel scheme for partitioning irregular graphs. {\em SIAM J. Sci. Comput.} {\bf 1998}, {\em 20}, 359--392.

\bibitem{Kimn:2007:ROB}
Kimn, J.-H.; Sarkis, M. Restricted overlapping balancing domain decomposition methods and restricted coarse problems for the {Helmholtz} problem. {\em Comput. Methods Appl. Mech. Engrg.} {\bf 2007}, {\em 196}, 1507--1514.

\bibitem{Kimn:2013:SLR}
Kimn, J.-H.; Sarkis, M. Shifted {Laplacian} {RAS} solvers for the {Helmholtz} equation. In {\em Domain Decomposition Methods in Science and Engineering XX}; Bank, R., Holst, M., Widlund, O., Xu, J., Eds.; Springer: Berlin, Germany, 2013; pp. 151--158.

\bibitem{Lahaye:2017:HTC}
Lahaye, D.; Vuik, C. How to choose the shift in the shifted {Laplace} preconditioner for the {Helmholtz} equation combined with deflation. In {\em Modern Solvers for Helmholtz Problems}; Lahaye, D., Tang, J., Vuik, K., Eds.; Birkh\"{a}user: Cham, Switzerland, 2017; pp. 85--112.

\bibitem{Lehoucq:1998:ARPACK}
Lehoucq, R.~B.; Sorensen, D.~C.; Yang, C. \textit{{ARPACK} Users' Guide: Solution of Large-Scale Eigenvalue Problems with Implicitly Restarted {A}rnoldi Methods}; Society for Industrial and Applied Mathematics (SIAM): Philadelphia, PA, USA, 1998.

\bibitem{Moiola:2014:SIREV}
Moiola, A.; Spence, E.~A. Is the Helmholtz equation really sign-indefinite? {\em SIAM Review} {\bf 2014}, {\em 56}, 274--312.

\bibitem{Nataf:2019:AGD}
Nataf, F.; Tournier, P.-H. A {GenEO} domain decomposition method for saddle point problems. \textit{arXiv preprint} {\bf 2019}, arXiv:1911.01858.

\bibitem{Nataf:2010:ATL}
Nataf, F.; Xiang, H.; Dolean, V. A two level domain decomposition preconditioner based on local {Dirichlet-to-Neumann} maps. {\em C. R. Math. Acad. Sci. Paris. I Math.} {\bf 2010}, {\em 348}, 1163--1167.

\bibitem{Nataf:2011:ACS}
Nataf, F.; Xiang, H.; Dolean, V.; Spillane, N. A coarse space construction based on local {Dirichlet-to-Neumann} maps. {\em SIAM J. Sci. Comput.} {\bf 2011}, {\em 33}, 1623--1642.

\bibitem{Spillane:2021:AAT}
Spillane, N. An abstract theory of domain decomposition methods with coarse spaces of the {GenEO} family. \textit{arXiv preprint} {\bf 2021}, arXiv:2104.00280.

\bibitem{Spillane:2014:ARC}
Spillane, N.; Dolean, V.; Hauret, P.; Nataf, F.; Pechstein, C.; Scheichl, R. Abstract robust coarse spaces for systems of {PDEs} via generalized eigenproblems in the overlaps. {\em Numer. Math.} {\bf 2014}, {\em 126}, 741--770.

\bibitem{Taus:2020:LSweeps}
Taus, M.; Zepeda-N{\'u}{\~n}ez, L.; Hewett, R.~J.; Demanet, L. {L-Sweeps}: A scalable, parallel preconditioner for the high-frequency {Helmholtz} equation. {\em J. Comput. Phys.} {\bf 2020}, {\em 420}, 109706.

\bibitem{Wang:2011:OMO}
Wang, S.; de Hoop, M.~V.; Xia, J. On {3D} modeling of seismic wave propagation via a structured parallel multifrontal direct {Helmholtz} solver. {\em Geophys. Prospect.} {\bf 2011}, {\em 59}, 857--873.

\bibitem{Zarmi:2013:AGA}
Zarmi, A.; Turkel, E. A general approach for high order absorbing boundary conditions for the {Helmholtz} equation. {\em J. Comput. Phys.} {\bf 2013}, {\em 242}, 387--404.

\end{thebibliography}
\end{document}